\documentclass[10pt]{article} 
\usepackage[utf8]{inputenc}
\usepackage[T1]{fontenc}
\usepackage{lmodern}
\usepackage{amsmath}
\usepackage{amssymb}
\usepackage{amsthm}
\usepackage{graphics}
\usepackage{graphicx}
\usepackage{stmaryrd}
\usepackage{xcolor,colortbl}
\definecolor{mycolor}{rgb}{.8, .8, 1}
\usepackage{caption}
\usepackage{placeins}
\usepackage{bbm}
\usepackage{multirow}
\usepackage{booktabs}
\usepackage{comment}
\usepackage{utfsym}
\usepackage{tikz}

\usepackage{subcaption}
\usepackage[ruled]{algorithm2e}
\usepackage{setspace} 

\usepackage{hyperref}
\hypersetup{
    colorlinks,
    citecolor=black,
    filecolor=black,
    linkcolor=black,
    urlcolor=black
}

\usepackage[top=1.5in, bottom=1.5in, left=1.5in, right=1.5in]{geometry}

\def\R{ { \mathbb{R} } }

\def\bbN{ { \mathbb{N} } }

\def\tref{\text{ref}}

\def \V{ { v } } 
\def \U{ { u } } 

\def \bfV{ { \mathbf{\V}  } }
\def \bfU{ { \mathbf{\U}  } }

\def \bfQ{ { \mathbf{q}  } }
\def \bfP{ { \mathbf{p}  } }
\def \bfq{ { \bar{\mathbf{q}}  } }
\def \bfp{ { \bar{\mathbf{p}}  } }

\def \bfY{ { \mathbf{y}  } }
\def \bfy{ { \overline{\mathbf{y}}  } }

\def \u{ { \overline{u} } }

\def \hami{ { \mathcal{H}  } }
\def \calD{ { \mathcal{D}  } }
\def \calE{ { \mathcal{E}  } }

\def \half{ { \tfrac{1}{2} } }

\def \AE{\text{AE}}

\def \model{\text{pred}}

\def \modelr{{\text{p}\overline{\text{re}}\text{d}}}

\def \stab{{\text{s}\overline{\text{ta}}\text{b}}}

\def \init{\text{init}}

\def \pred{\text{pred}}
\def \flux{\text{Flow}}

\newtheorem{remark}{Remark}

\newcommand{\sci}[1]{\!\times\! 10^{#1}}

\title{Hamiltonian reduction using a convolutional auto-encoder coupled to a  Hamiltonian neural network}

\usepackage{authblk}
\author[1]{Raphaël Côte}
\author[1, 2]{Emmanuel Franck}
\author[1, 2]{Laurent Navoret}
\author[1, 2]{Guillaume Steimer}
\author[1, 2]{Vincent Vigon}
\affil[1]{Institut de Recherche Mathématique Avancée, UMR 7501 Université de Strasbourg et CNRS, 7 rue René Descartes 67000 Strasbourg, France}
\affil[2]{INRIA, MACARON Project, Strasbourg, France}

\date{}

\newcommand{\answerG}[1]{\textcolor{black}{#1}} 

\begin{document}

\maketitle 


\begin{abstract}
    The reduction of Hamiltonian systems aims to build smaller reduced models, valid over a certain range of time and parameters, in order to reduce computing time. By maintaining the Hamiltonian structure in the reduced model, certain long-term stability properties can be preserved. In this paper, we propose a non-linear reduction method for models coming from the spatial discretization of partial differential equations: it is based on convolutional auto-encoders and Hamiltonian neural networks. Their training is coupled in order to simultaneously learn the encoder-decoder operators and the reduced dynamics. Several test cases on non-linear wave dynamics show that the method has better reduction properties than standard linear Hamiltonian reduction methods. 

\end{abstract}
\vspace{1cm}

\\


\clearpage

\section{Introduction}

Hamiltonian reduced order modeling techniques have been successfully developed in order to perform accelerated numerical simulations of some parameterized Hamiltonian models of large dimension \cite{peng2015symplectic, hesthaven2022reduced, ROmvlasov}. The spatial discretization of some wave-like partial differential equations gives rise to very large such Hamiltonian systems. Reduced order models can be essential for real-time simulations or when a large number of simulation instances are required as part of a control, optimisation or uncertainty quantification algorithm. Starting from the initial model, a large differential system, the methods consist into constructing a differential system of a smaller size that can produce valid approximate solutions for a predefined range of times and parameters. Many physical models have a Hamiltonian structure  and this gives the system a certain number of geometrical properties like the conservation of energy \answerG{ and the symplecticity of the phase space flows}. In particular, the preservation of this structure at the discrete level enables to ensure large-time stability of the numerical simulations \cite{Hairer}. In order to build consistent and robust reduced models, it is therefore interesting to preserve this Hamiltonian structure through the reduction. 
The construction of reduced models can be  divided in two steps: (i) find a so-called pair of encoder and decoder operators that goes from the full  to the reduced variables and inversely; (ii) identify the dynamics followed by the reduced variables. 
The construction of the encoder and decoder operators relies on a large number of data produced by numerical simulations in the range of time and parameters of interest. 

The first approach to reduce a large Hamiltonian system relies on a linear approximation: the solutions manifold is approximated with a \answerG{symplectic vector space} of small dimension \cite{peng2015symplectic}. The encoder is here a linear mapping, which is also \answerG{constructed to be symplectic} so that the Hamiltonian structure is preserved into the reduced model. Such symplectic mapping can be constructed from data through greedy algorithms \answerG{\cite{Afkham_2017}} or through a Singular Value Decomposition (SVD) methodology: this is the Proper Symplectic Decomposition (PSD) proposed in \cite{peng2015symplectic}. In this work, several algorithms have been proposed to define approximated optimal symplectic mappings: for instance, the cotangent-lift algorithm devise a symplectic mapping which is also orthogonal and have a block diagonal structure. Then the reduced model is obtained using the Galerkin projection method: the model is constructed by supposing that a symplectic projection of the residual vanishes, where the residual stands for the error obtained after replacing the original variables by the decoded reduced variables.

Such linear reductions, however, can hardly handle non-linear dynamics: this is the case for convection-dominated or non-linear wave like problems for which the solution manifold is badly approximated by hyperplanes. In order to build more expressive reduced models, one possibility is to consider time adaptive reduced methods \cite{pagliantini2021dynamical,hesthaven2023adaptive}. Another widely investigated possibility is to consider non-linear reduction methods. 

Regarding the construction of non-linear encoder-decoder operators, a first class of methods rely on manifold learning techniques \cite{sonday2010,hyperelastic2016}. Such methods are based on the geometrical analysis of the neighbors graph of the data thanks to the computation of geodesic distances (ISOMAP method, \cite{TenenbaumEtAl2000}), of eigenfunctions of the graph Laplacian (EigenMaps method \cite{belkin2003laplacian}) or of diffusion processes (DiffusionMaps method \cite{COIFMAN20065}). This provides reduced variables for each data that can be further interpolated using the  Nystr\"{o}m formula \cite{bengio2004learning}. 

Since the explosion of deep learning in the early 2010s, new dimension reduction methods grounded on neural networks have been developed. The convolutional auto-encoder architecture \cite{goodfellow2016deep} seems particularly appropriate since its very purpose is to determine latent variables: the neural network is indeed divided into an encoder part  and decoder part and they are trained simultaneously so that the sequence of encoder and decoder is close to the identity map. This was originally developed for image generation, but has been also used for reduced order modeling for models coming from the spatial discretization of partial differential equations on a grid \cite{AEmanifold, kim2022fast, romor2023non}. \answerG{Indeed, convolutional neural networks have proved particularly effective to  extract multi-scale informations of grid-structured data.  Secondly, they involve far less parameters than their dense counterparts, especially when the input size of the neural network is large.} In \cite{buchfink2021symplectic}, the authors use auto-encoder neural networks for Hamiltonian reduction: the encoder and the decoder are weakly constrained to be symplectic thanks to a penalization term in the cost functional. 

Once non-linear encoders and decoders have been devised from data, the dynamics of the reduced variables still has to be determined. Two strategies can be considered.  The first one relies on a Galerkin projection of the Hamiltonian system as in the linear case \cite{buchfink2021symplectic}. Note that the reduced model is indeed Hamiltonian provided the decoder is a symplectic map. However, the reduced model still requires the evaluation of the vector field in the original large dimension space of size $2N$: this is a well-known difficulty in non-linear reduction. To overcome this difficulty, hyper-reduction methods have been proposed like the discrete empirical interpolation method (DEIM) \cite{chaturantabut2010nonlinear}. \answerG{This method has to be adapted to not destroy the geometric structure of the full order model as in \cite{hesthaven2023adaptive, pagliantini2023fully,  Pagliantini_2023} where the authors propose to apply a DEIM algorithm to the gradient of the Hamiltonian thus, with some additional strategies, preserve the geometric structure of the full order model.}

Another approach is to learn about the dynamics of the reduced variables using a neural network: given the initial state, the neural network provides the full trajectory.  As the learning is done directly in the reduced dimension, the obtained reduced model does not require an evaluation of non-linear terms in the original variables: this is a clear advantage of the method compared with projection-based ones. The reduced dynamics can be captured for instance by Recurrent Neural Networks (RNNs), Long Short Term Memory (LSTM) neural networks \cite{maulik2021reduced}, or by fully connected networks  \cite{fresca2021comprehensive}. This has also been considered as correction of the Galerkin-type reduced models  \cite{wang2020recurrent}.

Here we consider another strategy which consists in learning the vector field that generates the observed reduced dynamics. The neural network is trained so as to minimize its deviation from the finite difference time derivative of the reduced data obtained after encoding.  As we aim at conserving the Hamiltonian structure at the reduced level, the vector field is further supposed to be associated with a reduced Hamiltonian function. Therefore, we can learn directly the Hamiltonian function instead of the vector field. This is a so-called Hamiltonian Neural Network (HNN) strategy proposed in \cite{hamiltonian_nn} where a symplectic time integrator is used. \answerG{We also note that   neural networks methods has also recently be used to learn hidden or reduced dynamics which also involve dissipation \cite{Yu_2021,Zhang_2022,chen2023constructing,Flaschel_2023}.}


The present paper proposes to combine an auto-encoder strategy for the encoding-decoding part and a HNN method to learn the reduced dynamics: this will be referred to as the AE-HNN method. Note that there is a priori no reason for the auto-encoder neural networks to spontaneously provide reduced variables compatible with Hamiltonian dynamics. Therefore, some constraints on the auto-encoder have to be added. This can be done by imposing symplecticity \answerG{weakly as in \cite{buchfink2021symplectic}, where a penalization term of the symplectic constraint is added to the loss function.} Here, we propose instead to train it simultaneously with the HNN.  With this joint training, the auto-encoder will gradually converge to a set of reduced variables compatible with a Hamiltonian system. 
Of course, this means that the loss functions associated with each neural network must be weighted judiciously during training. Note that such a joint training of the encoding-decoding operators and the reduce dynamics have been explored in \cite{fresca2021comprehensive,Yu_2021}, but without considering Hamiltonian structures for the first \answerG{and without symplectic time integrator for the second.} 


The outline of the article is as follows. In Section~\ref{section:param_hamiltonian_system}, we introduce parameterized Hamiltonian systems as well as the main steps for the construction of reduced order models. Section~\ref{section:non-linear-Hamiltonian-reduction} then presents the non-linear AE-HNN reduction method. In particular we describe the architectures and the loss functions used for the trainings. Finally, Section~\ref{section: results} is devoted to the numerical results: we apply our reduction method on the Hamiltonian systems obtained after spatial discretization of linear, non-linear wave equations and a shallow water system and compare it with the linear PSD reduction technique.

\section{Parameterized Hamiltonian systems and reduction}
\label{section:param_hamiltonian_system}

In this section, we introduce the notations used for the parameterized Hamiltonian systems and the main steps for the construction of a reduced model.

In the following, we often write vectors of interest with bold script letters, operators with capital italic letters, with their parameters as indices, and overline quantities when related to the reduced model.

\subsection{Parameterized Hamiltonian dynamics}

We consider a parameterized autonomous Hamiltonian system, whose solution, $\bfY(t;\mu) \in \R^{2N}$ with $N \in \bbN^\ast$, depends on time $t \in [0,T]$, with $T > 0$, and on a parameter $\mu \in \Xi  \subset \R^d$, with $d \in \mathbb{N}$. The dynamics derive from a given Hamiltonian function $\hami : \R^{2N} \times \Xi \to \R$ and writes
\begin{equation}
    \begin{cases}
    \displaystyle \frac{d}{dt} \mathbf{y} (t; \mu) = J_{2N} \, \nabla_{\mathbf{y}} \hami \left( \bfY(t;\mu); \mu \right), \quad \text{ in } (0,T], \\
    \bfY(0;\mu) = \bfY_{\init}(\mu), 
    \end{cases}\label{eq: symplectic formulation}
\end{equation}
where $\bfY_{\init}(\mu) \in \R^{2N}$ is a given initial condition and $J_{2N}$ refers to the canonical symplectic matrix 
\[ J_{2N} = \begin{pmatrix} 0_N & I_N \\ -I_N & 0_N \end{pmatrix}, \] 
with $I_N$ the identity matrix of dimension $N$. 

Introducing the canonical coordinates $\bfY = (\bfQ, \bfP)^T$, the system becomes:
\begin{align}
    \begin{cases}
    \displaystyle \frac{d}{dt} \bfQ (t; \mu) =  \nabla_{\bfP} \hami \left( \bfQ, \bfP ; \mu \right), & \text{ in } (0,T],\\[4pt] \displaystyle \frac{d}{dt} \bfP (t; \mu) = -\nabla_{\bfQ} \hami \left( \bfQ, \bfP ; \mu \right), & \text{ in } (0,T], \\
    \bfQ(0;\mu) = \bfQ_{\init}(\mu),&\\  \bfP(0;\mu) = \bfP_{\init}(\mu),&
    \end{cases}\label{eq: symplectic formulation PQ}
\end{align}
with $\bfQ_{\init}, \bfP_{\init} \in \R^N$ such that $\bfY_{\init} = (\bfQ_{\init}, \bfP_{\init})^T$. A key property of such systems is that the associated flow is symplectic, meaning that $\phi_t \left(\bfY_{\init}(\mu); \mu\right) = \bfY(t;\mu)$ satisfies the relation
\[ \left(\nabla_\bfY \phi_t \left(\bfY_{\init}(\mu); \mu\right)\right)^T J_{2N}\left(\nabla_\bfY \phi_t \left(\bfY_{\init}(\mu); \mu\right)\right) = J_{2N}. \]
One consequence is that the Hamiltonian $\hami$ is preserved along the flow
\[\forall t \in (0,T], \mu \in \Xi, \quad \hami \left( \bfY(t;\mu); \mu \right) = \hami \left( \bfY_{\init}(\mu); \mu \right), \]
which is of particular importance when considering physical systems.

In this work, we are specifically interested in Hamiltonian systems resulting from the space discretization of one-dimensional \answerG{and two-dimensional} wave-type equations. In such systems, $\bfQ \in \R^N$ refers to the height of the wave at grid points and $\bfP \in \R^N$ to the velocity of the wave also at grid points. Examples will be detailed in the numerical section.

In order to provide numerical approximations of the solution, specific numerical schemes have been developed to ensure the symplectic property at the discrete level \cite{Hairer}. These schemes also guarantee large time stability of the numerical solutions. Here we consider the standard second-order Störmer-Verlet scheme. Denoting $\bfY^n_\mu = \left( \bfQ^n_\mu, \bfP^n_\mu \right)^T \in \R^{2N}$ the approximate solution at time $t^n = n \Delta t$, with time step $\Delta t > 0$, one iteration of the scheme is defined by: 
\begin{equation}
    \begin{cases}
        \bfP^{n + \half}_\mu = \bfP^{n}_\mu - \half \Delta t \,\nabla_{\bfQ} \hami \left( \bfQ^{n}_\mu, \bfP^{n + \half}_\mu; \mu \right), \\[10pt]
        \bfQ^{n + 1}_\mu = \bfQ^{n+ \half}_{\mu} + \Delta t\, \left[ \nabla_{\bfP} \hami \left( \bfQ^{n}_\mu, \bfP^{n + \half}_\mu ; \mu \right) + \nabla_{\bfP} \hami \left( \bfQ^{n+1}_\mu, \bfP^{n + \half}_\mu ; \mu \right) \right], \\[10pt]
        \bfP^{n + 1}_\mu = \bfP^{n + \half}_\mu - \half \Delta t\, \nabla_{\bfQ} \hami \left( \bfQ^{n+1}_\mu, \bfP^{n + \half}_\mu; \mu \right).
    \end{cases}\label{eq: H leapfrog}
\end{equation}
Under the further assumption that the Hamiltonian $\hami$ is separable, i.e. is the sum of a function depending only $\bfQ$ and another depending only $\bfP$:
\[
    \hami(\bfY; \mu) = \hami^1(\bfQ; \mu) + \hami^2(\bfP; \mu),
\]
the implicit first two steps of \eqref{eq: H leapfrog} become explicit and the scheme simplifies into:
\begin{equation}
    \begin{cases}
        \bfP^{n + \half}_{\mu} = \bfP^{n}_{\mu} - \half \Delta t\, \nabla_{\bfQ} \hami^1 \left( \bfQ^{n}_{\mu}; \mu \right), \\[10pt]
        \bfQ^{n + 1}_{\mu} = \bfQ^{n}_{\mu} + \Delta t\, \nabla_{\bfP} \hami^2 \left(\bfP^{n + \half}_{\mu}; \mu \right), \\[10pt]
        \bfP^{n + 1}_{\mu} = \bfP^{n + \half}_{\mu} - \half \Delta t\, \nabla_{\bfQ} \hami^1 \left( \bfQ^{n+1}_{\mu}; \mu \right).
    \end{cases}\label{eq: H leapfrog explicit}
\end{equation}


\subsection{Hamiltonian reduced order modeling}\label{section: Hamiltonian reduced order modeling}

Solving Hamiltonian systems with large dimension $2N \gg 1$ numerically can be relatively costly, and this is especially true when we want to solve a large number of them for a parametric study, for example. Therefore, methods have been developed in order to construct reduced Hamiltonian systems of smaller size $2K \ll 2N$, which capture the main dynamics for a range of times $t$ and reduction parameters $\mu$.

We first have to define an appropriate change of variable. 
To do that, we search for a $2K$-dimensional manifold $\widehat{\mathcal{M}}$ that approximates well the manifold 
\[
    \mathcal{M} = \left\{ \bfY(t;\mu) \text{ with } t \in [0,T], \mu \in \Xi  \right\} \subset \R^{2N}
\]
formed by \answerG{the values taken by the solutions  of the differential equation \eqref{eq: symplectic formulation}. The manifold structure results from the Cauchy-Lipshitz (Picard-Lindhöf) theorem with parameters under some regularity assumptions of the Hamiltonian}. The manifold $\widehat{\mathcal{M}}$ is defined thanks to a so-called decoding operator $\calD_{\theta_d}: \R^{2K} \to \R^{2N}$:
\[
   \widehat{\mathcal{M}} =  \left\{ \calD_{\theta_d} \left( \bfy \right)\text{ with } \bfy \in \R^{2K} \right\} \subset \R^{2N}.
\]
We also consider a pseudo-inverse operator $\calE_{\theta_e}  :\R^{2N} \to \R^{2K}$, called the encoder, which satisfies the relation 
\[ \calE_{\theta_e} \circ \calD_{\theta_d} = \text{Id}_{\R^{2K}}. \]
To determine $\mathcal{D}_{\theta_d}$ and $\mathcal{E}_{\theta_e}$, we therefore ask for the projection operator $ \calD_{\theta_d} \circ \calE_{\theta_e}$ onto $\widehat{\mathcal{M}}$ to be close to the identity on a data set $U \subset \mathcal{M}$:
\begin{equation}
   \forall u \in U,\quad \calD_{\theta_d} \circ \calE_{\theta_e}(u) \approx u.
   \label{eq:identity}
\end{equation}
The data set $U$ is composed of snapshots of the solutions at different times  and various parameters, obtained with the symplectic algorithm defined above in \eqref{eq: H leapfrog}; it writes
\[
    U = \left\{
   \bfY^0_{\mu_1}, \dots ,\bfY^{M}_{\mu_1}, \dots ,\bfY^0_{\mu_P}, \dots, \bfY^{M}_{\mu_P}
    \right\},
\]
where $M \in \mathbb{N}^*$ is the number of time-step chosen and $P \in \mathbb{N}^*$ the number of sampled parameters.

\answerG{In addition to these approximation properties, we also ask for the reduced variables, 
\begin{equation*}
\bfy(t;\mu) = \calE_{\theta_e}(\bfY(t;\mu)) \in \R^{2K}, 
\end{equation*}
to follow a reduced Hamiltonian dynamics:}
\begin{equation}
    \begin{cases}
        \displaystyle \frac{d}{dt} \bfy(t;\mu) = J_{2K} \nabla_\bfy \overline{\hami}_{\theta_h} (\bfy(t;\mu); \mu), \quad \text{ in } (0,T], \\
        \bfy(0;\mu) = \calE_{\theta_e}(\bfY_{\init}(\mu)),
    \end{cases}\label{eq: reduced hamiltonian model}
\end{equation}
where $\overline{\hami}_{\theta_h} : \R^{2K} \times \Xi \to \R$ is a reduced Hamiltonian to be built. 

The most common approach for Hamiltonian reduced order modeling is called the Proper Symplectic Decomposition (PSD) \cite{peng2015symplectic}.  This method is briefly described in Appendix~\ref{section: PSD}. Although efficient for linear dynamics, it fails into reducing non-linear ones. This is why several non-linear Hamiltonian reduction techniques have been developed \cite{buchfink2021symplectic,hesthaven2023adaptive}. In the next section, we present a strategy based on the coupling of an Auto-Encoder (AE) and a Hamiltonian Neural Network (HNN) method.

\section{A non-linear Hamiltonian reduction method}
\label{section:non-linear-Hamiltonian-reduction}

The method proposed in this work consists in constructing the Hamiltonian reduced model via neural networks. More precisely, we aim at defining the following three neural networks:
\begin{itemize}
    \item a decoder $\mathcal D_{\theta_d}:\R^{2K} \to \R^{2N}$,
    \item an encoder $\mathcal E_{\theta_e}:\R^{2N} \to \R^{2K}$,
    \item a reduced Hamiltonian $\overline{\mathcal H}_{\theta_h}:\R^{2K}\times\Xi \to \R$, 
\end{itemize}
where $(\theta_d, \theta_e, \theta_h)$ stands for their parameters, such that the resulting reduced dynamics provides a good approximation of the initial one. An auto-encoder strategy will be used to define $\mathcal D_{\theta_d}$ and $\mathcal E_{\theta_e}$ while a Hamiltonian Neural Network will be considered for $\overline{\mathcal H}_{\theta_h}$. Note that the encoder and decoder are not enforced to be symplectic but the reduced model is. 

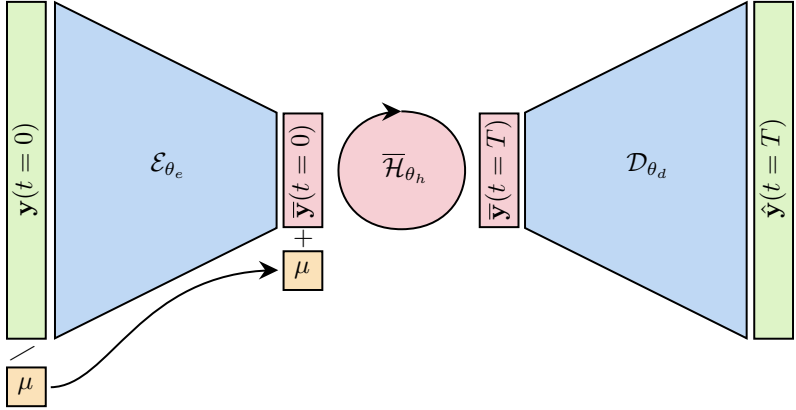
\begin{figure}[htb!]
\centering

\tikzset{every picture/.style={line width=0.75pt}}

\begin{tikzpicture}[x=0.75pt,y=0.75pt,yscale=-0.85,xscale=0.85]

\draw  [fill={rgb, 255:red, 184; green, 233; blue, 134 }  ,fill opacity=0.46 ] (1.03,0.67) -- (24.03,0.67) -- (24.03,200.68) -- (1.03,200.68) -- cycle ;
\draw  [fill={rgb, 255:red, 245; green, 166; blue, 35 }  ,fill opacity=0.32 ] (1.03,217.75) -- (24.03,217.75) -- (24.03,240.75) -- (1.03,240.75) -- cycle ;
\draw  [fill={rgb, 255:red, 74; green, 144; blue, 226 }  ,fill opacity=0.35 ] (29.5,0.93) -- (161,66.43) -- (161,134.93) -- (29.5,200.43) -- cycle ;
\draw  [fill={rgb, 255:red, 208; green, 2; blue, 27 }  ,fill opacity=0.19 ] (165.13,66.93) -- (188.13,66.93) -- (188.13,134.43) -- (165.13,134.43) -- cycle ;
\draw  [fill={rgb, 255:red, 245; green, 166; blue, 35 }  ,fill opacity=0.32 ] (165.13,148.75) -- (188.13,148.75) -- (188.13,171.75) -- (165.13,171.75) -- cycle ;
\draw  [fill={rgb, 255:red, 208; green, 2; blue, 27 }  ,fill opacity=0.19 ] (281.52,66.93) -- (304.52,66.93) -- (304.52,134.43) -- (281.52,134.43) -- cycle ;
\draw  [fill={rgb, 255:red, 74; green, 144; blue, 226 }  ,fill opacity=0.35 ] (439.9,200.43) -- (308.4,134.93) -- (308.4,66.43) -- (439.9,0.92) -- cycle ;
\draw  [fill={rgb, 255:red, 184; green, 233; blue, 134 }  ,fill opacity=0.46 ] (444.4,0.67) -- (467.4,0.67) -- (467.4,200.68) -- (444.4,200.68) -- cycle ;
\draw [fill={rgb, 255:red, 208; green, 2; blue, 27 }  ,fill opacity=0.19 ]   (235,65.67) .. controls (285.4,66.07) and (285.4,135.27) .. (235.4,135.67) .. controls (186.4,136.07) and (185.04,69.23) .. (232.06,65.8) ;
\draw [shift={(235,65.67)}, rotate = 179.08] [fill={rgb, 255:red, 0; green, 0; blue, 0 }  ][line width=0.08]  [draw opacity=0] (11.61,-5.58) -- (0,0) -- (11.61,5.58) -- (7.71,0) -- cycle    ;
\draw [line width=0.75]    (27,229.5) .. controls (62.15,229.01) and (79.65,161.37) .. (159.56,160.02) ;
\draw [shift={(162,160)}, rotate = 180] [fill={rgb, 255:red, 0; green, 0; blue, 0 }  ][line width=0.08]  [draw opacity=0] (11.61,-5.58) -- (0,0) -- (11.61,5.58) -- (7.71,0) -- cycle    ;

\draw (2.43,129.68) node [anchor=north west][inner sep=0.75pt]  [rotate=-270]  {$\mathbf{y}( t=0)$};
\draw (5.53,222.15) node [anchor=north west][inner sep=0.75pt]    {$\mu$};
\draw (166.52,130.68) node [anchor=north west][inner sep=0.75pt]  [rotate=-270]  {$\overline{\mathbf{y}}( t=0)$};
\draw (169.63,153.15) node [anchor=north west][inner sep=0.75pt]    {$\mu $};
\draw (169.13,134.4) node [anchor=north west][inner sep=0.75pt]    {$+$};
\draw (12.05,198.18) node [anchor=north west][inner sep=0.75pt]  [rotate=-40]  {$/$};
\draw (282.92,131.68) node [anchor=north west][inner sep=0.75pt]  [rotate=-270]  {$\overline{\mathbf{y}}( t=T)$};
\draw (445.8,130.68) node [anchor=north west][inner sep=0.75pt]  [rotate=-270]  {$\hat{\mathbf{y}}( t=T)$};
\draw (84.4,88.58) node [anchor=north west][inner sep=0.75pt]  [font=\normalsize]  {$\mathcal{E}_{\theta _{e}}$};
\draw (366.8,88.58) node [anchor=north west][inner sep=0.75pt]  [font=\normalsize]  {$\mathcal{D}_{\theta _{d}}$};
\draw (222,88.57) node [anchor=north west][inner sep=0.75pt]  [font=\normalsize]  {$\overline{\mathcal{H}}_{\theta _{h}}$};

\end{tikzpicture}

\caption{Prediction using the reduced model. The closed loop in the middle refers to the application of several iterations of the Störmer-Verlet scheme.}\label{fig: blueprint}
\end{figure}

Figure~\ref{fig: blueprint} illustrates how the reduced model is expected to be used for prediction. The initial condition $\bfY(t=0;\mu)$ is converted by the encoder $\calE_{\theta_e}$ to the reduced initial condition $\bfy(t=0;\mu)$. Then several iterations of the Störmer-Verlet scheme with the reduced Hamiltonian $\overline{\mathcal H}_{\theta_h}$ are performed to obtain an approximated reduced solution $\bfy(t=T;\mu)$ at time $T$. Finally, by using the decoder $\calD_{\theta_d}$, the latter is transformed into $\widehat \bfY(t=T;\mu) \approx \bfY(t=T;\mu)$. Note parameter $\mu$ has to be supplied to the Hamiltonian function. 

In order to determine the appropriate parameters of the three neural networks $\mathcal D_{\theta_d}$,  $\mathcal E_{\theta_e}$ and  $\overline{\mathcal H}_{\theta_h}$, we have to define both their architectures and the loss functions used for their training. This section focuses on the latter.

\FloatBarrier

\subsection{Reduction with an Auto Encoder (AE)}

An auto-encoder (AE) is a classical architecture of neural networks to find a reduced representation of data \cite{goodfellow2016deep}. It is composed of two neural networks, $\mathcal D_{\theta_d}$ and  $\mathcal E_{\theta_e}$, which are trained together such as to make the projection operator $\mathcal{D}_{\theta_d} \circ \mathcal{E}_{\theta_e}$ the closest to the identity map on the training data set $U$. Therefore, the AE is trained so as to minimize the following loss
\begin{equation}
    \mathcal{L}_{\AE} (\theta_e, \theta_d) = \sum_{\bfY \in \,U} \left\| \bfY -  \mathcal{D}_{\theta_d} \left( \mathcal{E}_{\theta_e} \left( \bfY \right) \right) \right\|^2_2.
    \label{eq:AE_loss}
\end{equation}
To account for the particular structure of $\bfY$ made of coordinates and momenta, the encoder input is a tensor of size $(N,2)$. This AE will be referred to as the bichannel AE.

Another choice would be to define two separate auto-encoders for coordinates and momenta: the coordinates AE is denoted $(\mathcal{E}^{1}_{\theta_e,1}, \mathcal{D}^{1}_{\theta_d,1})$ and the momenta AE is denoted $(\mathcal{E}^{2}_{\theta_e,2}, \mathcal{D}^{2}_{\theta_d,2})$; each encoder input has shape $(N,1)$. The AEs are trained by minimizing the loss:
\[
    \mathcal{L}_{\text{split},\AE} (\theta_e, \theta_d) = \sum_{(\bfQ, \bfP) \in \,U} \left\| \bfQ -  \mathcal{D}_{\theta_d,1}^{1} \left( \mathcal{E}_{\theta_e,1}^{1} \left( \bfQ \right) \right) \right\|^2_2 + \left\| \bfP -  \mathcal{D}_{\theta_d,2}^{2} \left( \mathcal{E}_{\theta_e,2}^{2} \left( \bfP \right) \right) \right\|^2_2.
\]
These AEs will be called the split AE. \answerG{This split AE will be used to preserve the separability property of the Hamiltonian, where applicable (see Remark~\ref{rem:separability} below).}

\begin{figure}[htb!]
\centering

\tikzset{every picture/.style={line width=0.75pt}} 

\begin{tikzpicture}[x=0.75pt,y=0.75pt,yscale=-1,xscale=1]

\draw  [fill={rgb, 255:red, 226; green, 226; blue, 226 }  ,fill opacity=1 ] (13,24) -- (33.67,24) -- (33.67,171.57) -- (13,171.57) -- cycle ;
\draw  [color={rgb, 255:red, 0; green, 75; blue, 170 }  ,draw opacity=1 ][fill={rgb, 255:red, 226; green, 226; blue, 226 }  ,fill opacity=1 ][line width=1.5]  (63.41,24) -- (83.92,24) -- (83.92,171.57) -- (63.41,171.57) -- cycle ;
\draw  [color={rgb, 255:red, 0; green, 75; blue, 170 }  ,draw opacity=1 ][fill={rgb, 255:red, 226; green, 226; blue, 226 }  ,fill opacity=1 ][line width=1.5]  (53,201.89) -- (94.33,201.89) -- (94.33,273.68) -- (53,273.68) -- cycle ;
\draw  [fill={rgb, 255:red, 226; green, 226; blue, 226 }  ,fill opacity=1 ] (255.8,215.27) -- (274.33,215.27) -- (274.33,261.64) -- (255.8,261.64) -- cycle ;
\draw  [color={rgb, 255:red, 0; green, 75; blue, 170 }  ,draw opacity=1 ][fill={rgb, 255:red, 226; green, 226; blue, 226 }  ,fill opacity=1 ][line width=1.5]  (123.67,202.56) -- (164.33,202.56) -- (164.33,274.35) -- (123.67,274.35) -- cycle ;
\draw  [color={rgb, 255:red, 0; green, 75; blue, 170 }  ,draw opacity=1 ][fill={rgb, 255:red, 226; green, 226; blue, 226 }  ,fill opacity=1 ][line width=1.5]  (195.16,165.33) -- (215.67,165.33) -- (215.67,312.91) -- (195.16,312.91) -- cycle ;
\draw  [color={rgb, 255:red, 110; green, 185; blue, 27 }  ,draw opacity=1 ][fill={rgb, 255:red, 226; green, 226; blue, 226 }  ,fill opacity=1 ][line width=1.5]  (314.49,164.67) -- (335,164.67) -- (335,312.24) -- (314.49,312.24) -- cycle ;
\draw  [color={rgb, 255:red, 110; green, 185; blue, 27 }  ,draw opacity=1 ][fill={rgb, 255:red, 226; green, 226; blue, 226 }  ,fill opacity=1 ][line width=1.5]  (365.67,202.56) -- (407,202.56) -- (407,274.35) -- (365.67,274.35) -- cycle ;
\draw  [color={rgb, 255:red, 110; green, 185; blue, 27 }  ,draw opacity=1 ][fill={rgb, 255:red, 226; green, 226; blue, 226 }  ,fill opacity=1 ][line width=1.5]  (436.33,202.56) -- (477,202.56) -- (477,274.35) -- (436.33,274.35) -- cycle ;
\draw  [color={rgb, 255:red, 110; green, 185; blue, 27 }  ,draw opacity=1 ][fill={rgb, 255:red, 226; green, 226; blue, 226 }  ,fill opacity=1 ][line width=1.5]  (446.33,24.67) -- (467,24.67) -- (467,172.24) -- (446.33,172.24) -- cycle ;
\draw  [fill={rgb, 255:red, 226; green, 226; blue, 226 }  ,fill opacity=1 ] (495.83,24.67) -- (516.33,24.67) -- (516.33,172.24) -- (495.83,172.24) -- cycle ;
\draw [color={rgb, 255:red, 110; green, 185; blue, 27 }  ,draw opacity=1 ][line width=1.5]    (314.49,164.67) -- (274.33,215.27) ;
\draw [color={rgb, 255:red, 0; green, 75; blue, 170 }  ,draw opacity=1 ][line width=1.5]    (215.67,164.67) -- (255.8,215.27) ;
\draw [color={rgb, 255:red, 0; green, 75; blue, 170 }  ,draw opacity=1 ][line width=1.5]    (215.67,312.24) -- (255.8,261.64) ;
\draw [color={rgb, 255:red, 110; green, 185; blue, 27 }  ,draw opacity=1 ][line width=1.5]    (274.33,261.64) -- (314.49,312.24) ;
\draw [color={rgb, 255:red, 0; green, 75; blue, 170 }  ,draw opacity=1 ][line width=1.5]    (164.33,202.56) -- (195.16,164.67) ;
\draw [color={rgb, 255:red, 0; green, 75; blue, 170 }  ,draw opacity=1 ][line width=1.5]    (164.33,274.35) -- (195.16,312.24) ;
\draw [color={rgb, 255:red, 110; green, 185; blue, 27 }  ,draw opacity=1 ][line width=1.5]    (335,164.67) -- (365.67,202.56) ;
\draw [color={rgb, 255:red, 110; green, 185; blue, 27 }  ,draw opacity=1 ][line width=1.5]    (365.67,274.35) -- (335,312.24) ;
\draw  [color={rgb, 255:red, 0; green, 75; blue, 170 }  ,draw opacity=0.3 ][line width=1.5]  (39,13) -- (251.5,13) -- (251.5,325) -- (39,325) -- cycle ;
\draw  [color={rgb, 255:red, 110; green, 185; blue, 27 }  ,draw opacity=0.3 ][line width=1.5]  (280.5,13.5) -- (493,13.5) -- (493,325.5) -- (280.5,325.5) -- cycle ;
\draw  [color={rgb, 255:red, 0; green, 0; blue, 0 }  ,draw opacity=0.4 ][dash pattern={on 5.63pt off 4.5pt}][line width=1.5]  (43.5,18.5) -- (99.5,18.5) -- (99.5,280) -- (43.5,280) -- cycle ;
\draw  [color={rgb, 255:red, 0; green, 0; blue, 0 }  ,draw opacity=0.4 ][dash pattern={on 5.63pt off 4.5pt}][line width=1.5]  (433.5,18) -- (489.5,18) -- (489.5,279.5) -- (433.5,279.5) -- cycle ;

\draw (16.71,174.71) node [anchor=north west][inner sep=0.75pt]  [font=\normalsize]  {$\mathbf{y}$};
\draw (258,265) node [anchor=north west][inner sep=0.75pt]  [font=\normalsize]  {$\overline{\mathbf{y}}$};
\draw (501.33,174.67) node [anchor=north west][inner sep=0.75pt]  [font=\normalsize]  {$\hat{\mathbf{y}}$};
\draw (16.33,119) node [anchor=north west][inner sep=0.75pt]  [font=\footnotesize,rotate=-270]  {$( N,2)$};
\draw (67,119) node [anchor=north west][inner sep=0.75pt]  [font=\footnotesize,rotate=-270]  {$( N,2)$};
\draw (449.67,119) node [anchor=north west][inner sep=0.75pt]  [font=\footnotesize,rotate=-270]  {$( N,2)$};
\draw (500.33,119) node [anchor=north west][inner sep=0.75pt]  [font=\footnotesize,rotate=-270]  {$( N,2)$};
\draw (64.83,260.5) node [anchor=north west][inner sep=0.75pt]  [font=\footnotesize,rotate=-270]  {$\left(\half N,4\right)$};
\draw (134.17,260.5) node [anchor=north west][inner sep=0.75pt]  [font=\footnotesize,rotate=-270]  {$\left(\half N,4\right)$};
\draw (376.17,260.5) node [anchor=north west][inner sep=0.75pt]  [font=\footnotesize,rotate=-270]  {$\left(\half N,4\right)$};
\draw (448.17,260.5) node [anchor=north west][inner sep=0.75pt]  [font=\footnotesize,rotate=-270]  {$\left(\half N,4\right)$};
\draw (198.33,255) node [anchor=north west][inner sep=0.75pt]  [font=\footnotesize,rotate=-270]  {$( 2N,)$};
\draw (317.67,255) node [anchor=north west][inner sep=0.75pt]  [font=\footnotesize,rotate=-270]  {$( 2N,)$};
\draw (258.33,254.5) node [anchor=north west][inner sep=0.75pt]  [font=\footnotesize,rotate=-270]  {$( 2K,)$};
\draw (171.33,262.33) node [anchor=north west][inner sep=0.75pt]  [rotate=-270] [align=left] {\text{flatten}};
\draw (340.67,270.83) node [anchor=north west][inner sep=0.75pt]  [rotate=-270] [align=left] {\text{unflatten}};
\draw (223.33,277.33) node [anchor=north west][inner sep=0.75pt]  [rotate=-270] [align=left] {\text{dense layers}};
\draw (288.67,277.33) node [anchor=north west][inner sep=0.75pt]  [rotate=-270] [align=left] {\text{dense layers}};
\draw (39.67,94.33) node [anchor=north west][inner sep=0.75pt]  [font=\LARGE]  {$\ast $};
\draw (100.67,230.67) node [anchor=north west][inner sep=0.75pt]  [font=\LARGE]  {$\ast $};
\draw (414,230.67) node [anchor=north west][inner sep=0.75pt]  [font=\LARGE]  {$\ast $};
\draw (473.67,94.33) node [anchor=north west][inner sep=0.75pt]  [font=\LARGE]  {$\ast $};
\draw (450,178) node [anchor=north west][inner sep=0.75pt]  [font=\large]  {$\Uparrow $};
\draw (66.67,178) node [anchor=north west][inner sep=0.75pt]  [font=\large]  {$\Downarrow $};
\draw (177.5,16.07) node [anchor=north west][inner sep=0.75pt]  [font=\normalsize,color={rgb, 255:red, 0; green, 75; blue, 170 }  ,opacity=1 ]  {$\text{encoder} \ \mathcal{E}_{\theta _{e}}$};
\draw (283.5,16.9) node [anchor=north west][inner sep=0.75pt]  [color={rgb, 255:red, 91; green, 148; blue, 35 }  ,opacity=1 ]  {$\text{decoder} \ \mathcal{D}_{\theta _{d}}$};
\draw (102,155.25) node [anchor=north west][inner sep=0.75pt]  [color={rgb, 255:red, 0; green, 0; blue, 0 }  ,opacity=0.6 ,rotate=-270] [align=left] {encoder block};
\draw (415,155.25) node [anchor=north west][inner sep=0.75pt]  [color={rgb, 255:red, 0; green, 0; blue, 0 }  ,opacity=0.6 ,rotate=-270] [align=left] {decoder block};

\end{tikzpicture}

\caption{Auto-encoder architecture: encoder in blue, decoder in green. Symbols. $\ast$: stride 1 convolution with periodic padding, $\Downarrow$ down-sampling (stride 2 convolution), $\Uparrow$: up-sampling (repeat once each value along the last axis then smooth it with a kernel size $2$ convolution).}

\label{fig: schema autoencoder}
\end{figure}
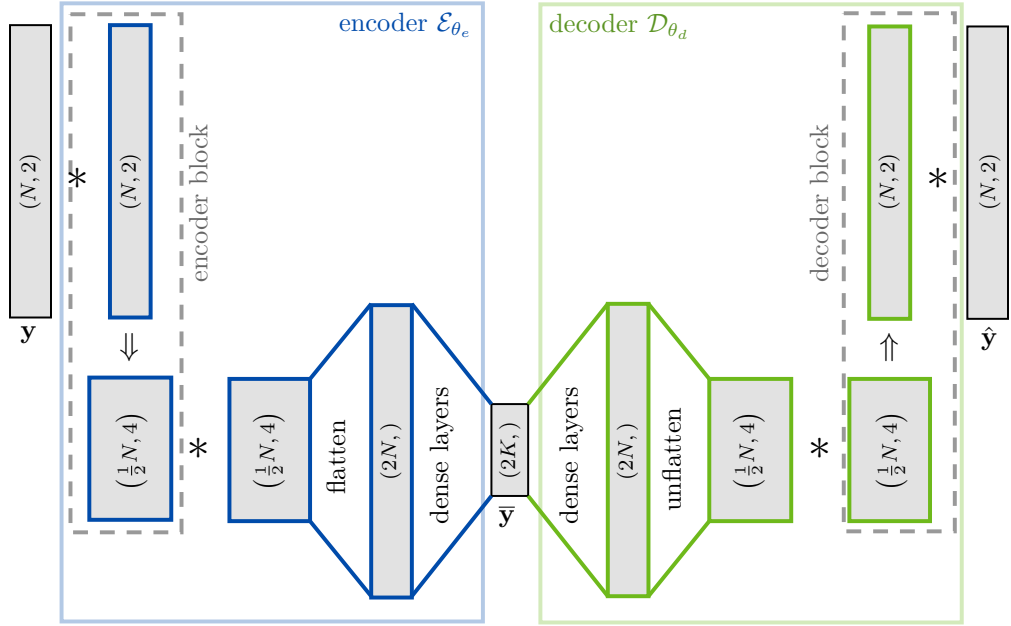

The architectures of the neural networks are chosen specifically to the Hamiltonian systems in consideration. In this work, we focus on systems resulting from the spatial discretization of wave-like equations: networks will be more efficient if they take into account the spatial structure of the data. Consequently, the encoder $\mathcal E_{\theta_e}$ is first composed of several pairs of convolution layer and down-sampling before ending with some dense layers, as depicted in Figure~\ref{fig: schema autoencoder}. As usual for AE networks, the decoder is constructed in a mirror way, i.e. starting with some dense layers and then ending with pairs of up-sampling and convolution layers  in reversed size order. 


More precisely, the encoder takes as input a vector $\bfY = (\bfQ, \bfP)$ and  starts with a succession of so-called encoder blocks, made of a stride $1$ convolution with kernel size $3$ and a down-sampling step (stride $2$ convolution with kernel size $2$). An encoder block results in an output that has twice as many channels and half as many rows as the input. After possibly composing several encoder blocks, we add a last convolution layer with kernel size $3$ and then a flattening operation by concatenating every channels. Then, dense layers are added until reaching the desired reduced dimension $2K$ of $\bfy$. As already said, the decoder is built as a mirror: the flattening operation is replaced by unflattening and the encoder blocks by the decoder blocks made of an up-sampling layer of size $2$ smoothed out with a convolution with kernel size $2$ and a convolution layer with kernel size $3$, which symmetrically results in output that has half as many channels and twice as many rows as the input. The architecture of the auto-encoder is thus defined with the number of encoder and decoder blocks and the dense layer sizes for both encoder and decoder. Figure~\ref{fig: schema autoencoder} illustrates an example of auto-encoder architecture with one block for the encoder and one block for the decoder.

\answerG{
This AE architecture can easily be extended to 2D systems using two-dimensional convolutional layers and up and down-sampling with appropriate dimensions.
}


\subsection{Reduced model with a Hamiltonian Neural Network (HNN)}\label{section: Reduced model with a Hamiltonian Neural Network}

The AE constructed in the previous section enables us to define the reduced trajectories:
\begin{equation}
    \bfy(t;\mu) = \mathcal{E}_{\theta_e}(\bfY(t;\mu)).
    \label{eq:reduced_trajectories}
\end{equation}
To obtain the dynamics of these reduced variables, we propose to use a Hamiltonian Neural Network strategy \cite{hamiltonian_nn}. We thus look for a neural network function $\overline\hami_{\theta_h}$, parameterized by $\theta_h$, such that: 
\[
    \frac{d}{dt} \bfy(t;\mu) = J_{2K} \nabla_{\bfy}\overline\hami_{\theta_h}(\bfy(t;\mu);\mu).
\]
Note that the Hamiltonian is supposed to depend on parameter $\mu$. We remind that $\mu \in \Xi$ stands for known parameters of the model, unlike $\theta_h$ that is the neural networks parameters to be learnt.
The architecture of the reduced Hamiltonian is a classical MLP neural network. The size of the neural network is chosen to be small so that the reduced model remains competitive. This reduced dynamics are in practice defined through a time discretization. We therefore introduce the  prediction operator: 
\begin{equation*}
    \mathcal{P}_s\left( \bfy; \overline\hami_{\theta_h,\mu} \right),
\end{equation*}
which consists in performing $s \in \mathbb{N}^*$ iterations of the Störmer-Verlet scheme, defined in \eqref{eq: H leapfrog}, starting from $\bfy$ and where $\overline\hami_{\theta_h,\mu}$ stands for the Hamiltonian function $\overline\hami_{\theta_h}(.;\mu)$. The number of steps $s$ considered in this prediction is called the watch duration. This is a hyper-parameter of the method that has to be set. The parameters of the reduced Hamiltonian are finally obtained by minimizing the following loss:
\begin{equation}
  \mathcal{L}^s_{\modelr} (\theta_e,\theta_h) = \sum_{ \bfY^{n}_{\mu} ,  \bfY^{n+s}_{\mu}  \in \, U}\left\| \bfy^{n+s}_{\mu} - \mathcal{P}_s\left( \bfy^{n}_{\mu} ; \overline\hami_{\theta_h,\mu} \right) \right\|^2_2
 \label{eq:HNN_loss}
\end{equation}
where $\bfY^{n}_{\mu} , \bfY^{n+s}_{\mu} \in \, U$ denote the sampling of random pairs $(\bfY^{n}_{\mu}, \bfY^{n+s}_{\mu})$ on the data set $U$. In other words, random time series of size $s$ are sampled and only the data at both ends are considered. This loss thus compares the reduced trajectories \eqref{eq:reduced_trajectories} with the ones obtained with the reduced Hamiltonian. The name of the loss function, ``$\modelr$'', refers to the prediction in the reduced variables. Note that the encoder neural network  is required to obtain the reduced data $\bfy = \mathcal{E}_{\theta_e}(\bfY)$ and this is why the loss also depends on $\theta_e$. This kind of loss function, based on a model, has been widely used in physics based deep learning methods  \cite{thuerey2021physics}.

Then, we constrain the reduced trajectories to preserve the reduced Hamiltonian with the following loss function : 
\begin{equation}
    \mathcal{L}^s_{\stab} (\theta_e, \theta_h) = \sum_{\bfY^{n}_{\mu}, \bfY^{n+s}_{\mu}\in \, U} \left\| \overline\hami_{\theta_h,\mu} \left( \bfy^{n+s}_{\mu} \right) - \overline\hami_{\theta_h,\mu} \left( \bfy^{n}_{\mu} \right) \right\|^2_2,
    \label{eq:energy_loss}
\end{equation}
where $\bfy^{n+s}_{\mu}$ and $\bfy^{n}_{\mu}$ are still obtained using the encoder $\mathcal E_{\theta_e}$. The aim is to ensure some stability of the reduced model, hence its name ``$\stab$'' . At first sight, this loss seems redundant with the prediction-reduced loss since using Störmer-Verlet schemes in the prediction step ensures that the reduced Hamiltonian is preserved at least approximately. Hence if the prediction-reduced loss becomes small, so does the stability loss. However, this additional loss may help the coupling to converge.

\begin{remark}\label{rem:separability}
\answerG{The separability of the Hamiltonian could be an interesting property to preserve at the reduced level. Using a split AE to learn separately reduced coordinates and momenta, a separable reduced Hamiltonian can be designed: $$\overline{\hami}_{\theta_h}(\bfy;\mu) = \overline{\hami}_{\theta_h}^1(\mathbf{\bar q};\mu) + \overline{\hami}_{\theta_h}^2(\mathbf{\bar p};\mu),$$ 
involving two neural networks. This will be referred to as the split HNN.  The Störmer-Verlet scheme would then have a cost of an explicit scheme. Note however that this is not a crucial gain since the reduced models under consideration have small sizes.}
\end{remark}

\begin{remark}
\answerG{The HNN will produce an approximation of the Hamiltonian, whose associated flow would have generated the discrete solution we would like to fit. Note also that the data used to fit the HNN was obtained after encoding a discrete solution of the initial dynamics of the Hamiltonian. Thus, even this discrete dynamics is an approximation to the encoded continuous Hamiltonian flow. However, from a practical point of view, it is not essential to capture the underlying continuous dynamics (at the reduced level): the objective is rather to obtain a method capable of reproducing the discrete dynamics, with good geometric properties, for a given time step. So, in practice, we actually do not vary the time step. }
\end{remark}

\subsection{Strong coupling of the neural networks} 

The prediction-reduced and the stability-reduced losses \eqref{eq:HNN_loss} already introduce a coupling between the encoder neural networks and the reduced Hamiltonian one. To make the coupling stronger, 
we could ask for the trajectories in the initial variables to be well recovered. This is why we introduce the following fourth loss function named ``\model'' which refers to the model prediction in the FOM space: 
\begin{equation}
\mathcal{L}_{\model}^s (\theta_e,\theta_d,\theta_h) = \sum_{\bfY^n_{\mu}, \bfY^{n+s}_{\mu} \in \, U} \left\| \bfY^{n+s}_{\mu} - \mathcal{D}_{\theta_d} \left( \mathcal{P}_s\left( \bfy^n_{\mu}; \overline\hami_{\theta_h,\mu} \right) \right) \right\|^2_2.
\label{eq:modelr}
\end{equation}
This is the only loss function that couples the three neural networks. It compares the trajectories in the initial variables with the full process of encoding, predicting in reduced variables over $s$ iterations and then decoding.

To sum up, we use four different loss functions $\mathcal{L}_{\AE}, \mathcal{L}_{\modelr}^s, \mathcal{L}_{\stab}$ and $\mathcal{L}_{\model}^s$, given by \eqref{eq:AE_loss}-\eqref{eq:HNN_loss}-\eqref{eq:energy_loss}-\eqref{eq:modelr} that couple both AE and HNN neural networks. The training aims to find the parameters $(\theta_e, \theta_d, \theta_h)$ that are a solution to the minimization problem:
\begin{equation*}
    \underset{\theta_e, \theta_d, \theta_h}{\operatorname{min}}\quad\omega_{\AE}\, \mathcal{L}_{\AE} (\theta_e,\theta_d) + \omega_{\modelr}\, \mathcal{L}^s_{\modelr}(\theta_e,\theta_h) + \omega_{\stab}\, \mathcal{L}^s_{\stab} (\theta_e, \theta_h) + \omega_{\model}\, \mathcal{L}^s_{\model} (\theta_e,\theta_d,\theta_h),
\end{equation*}
where $\omega_{\AE}$, $\omega_{\modelr}$, $\omega_{\stab}$, $\omega_{\model}$ are positive weights: these are hyper-parameters of the method. The four loss functions interact during training and possibly compete with each other. Note that in the end, the only loss value that really quantifies the quality of the reduction and prediction process is the one corresponding to $\mathcal{L}_{\model}^s$. The other loss functions are only useful in the training process.

\subsection{Training hyper-parameters}\label{section: learning details}

In addition to the parameters of the neural networks (number of layers, size of the layers), the training of the model also depends on several hyper-parameters. 


\answerG{\paragraph{Reduced dimension $K$.} In classical reduction method, the larger $K$, the more accurate the reduced model. Regarding the AE-HNN method, as the approximation is truly non-linear, there may be no benefit increasing the reduced dimension. In practice, the minimum possible reduced dimension should be equal to the number of variable parameters in the model.
}

\paragraph{Watch duration in predictions.} One of the hyper-parameter to set is the watch duration $s$ in the loss functions $\mathcal{L}^s_{\stab}$,  $\mathcal{L}^s_{\modelr}$ and $\mathcal{L}^s_{\model}$ that make predictions. This quantity should be not too small to capture the dynamics but also not too large as the computation of the gradients of the associated loss functions may generate vanishing gradient problems. In the numerical setting, the watch duration will be typically set to $ s=16$.

\begin{figure}[htb!]
    \begin{subfigure}{\textwidth}
        \centering
        \includegraphics[width=1.\textwidth]{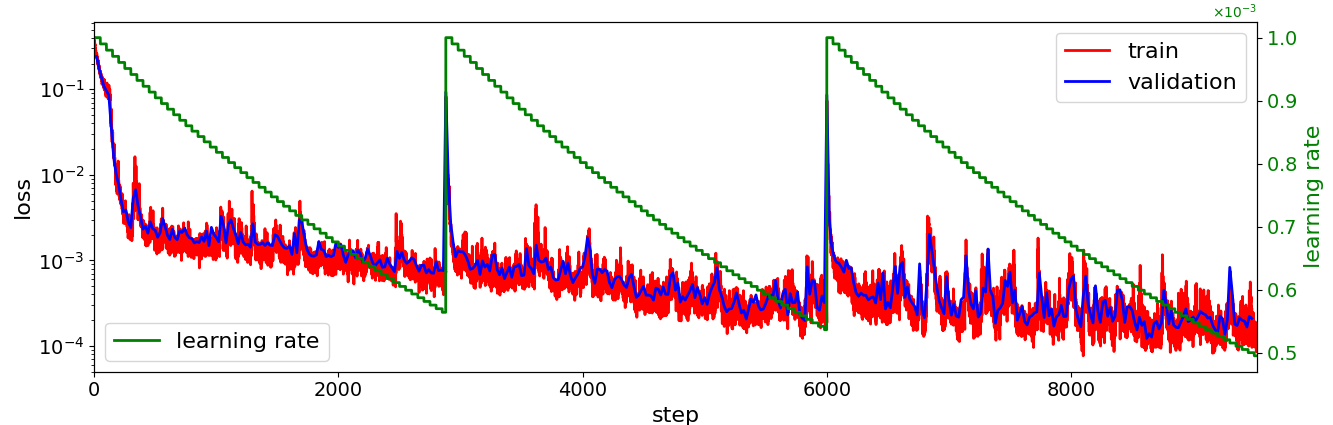}
        \caption{Training loss function (blue) and validation loss function (red) as functions of the training step.}
        \label{fig: training losses example 1}
    \end{subfigure}
    \begin{subfigure}{\textwidth}
        \centering
        \includegraphics[width=1.\textwidth]{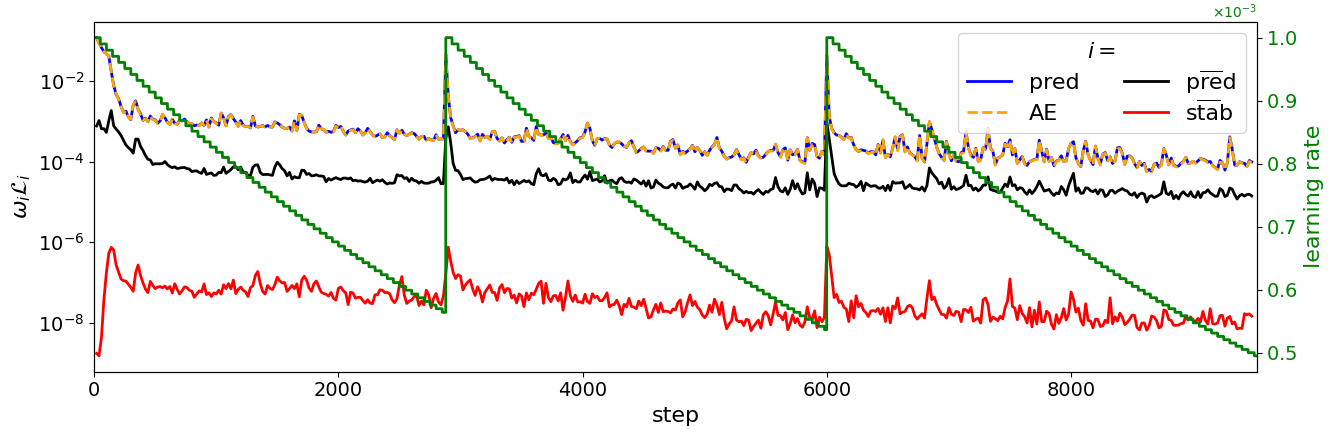}
        \caption{All the weighted loss functions as functions of the training step.}
        \label{fig: training losses example 2}
    \end{subfigure}
            
    \caption{Example of loss functions history during a training, overlaid with the evolution of the learning rate (green).}
    \label{fig: training losses example}
\end{figure}

\paragraph{Loss functions weights.} Losses weights have been chosen experimentally as follows:
\[ \omega_{\model}= 0.1,\quad \omega_{\AE}=0.1,\quad \omega_{\modelr}=80,\quad \omega_{\stab}=7\times 10^{-4}. \]
A typical loss functions history is shown in Figure \ref{fig: training losses example 2}: each loss function is represented multiplied by its weight. We first notice that the prediction loss function $\mathcal{L}^s_{\model}$ and the auto-encoder loss function $\mathcal{L}_{\AE}$ have the same magnitude and actually are almost equal: the error in prediction is mostly due to the encoder-decoder step. We keep this behavior by assigning them the same weight $\omega_{\model}=\omega_{\AE}=0.1$. This value is determined in proportion to the learning rate chosen below. As the prediction loss function $\mathcal{L}^s_{\model}$ is the most important one for the applications, we want it to dominate over the others. We therefore set the weight $\omega_{\modelr}=80$ so that the weighted reduced prediction loss function $\omega_{\modelr}\mathcal{L}^s_{\modelr}$ is about $10$ times smaller than the previous two. Finally, we want the weighted reduced stability loss function $\omega_{\stab}\mathcal{L}^s_{\stab}$ to act as a quality control that remains small compared to the other loss functions. To this end, we set $\omega_{\stab}=7\times 10^{-4}$ in order to make it about $100$ times smaller than the weighted reduced prediction loss.


\paragraph{Gradient descent.} An Adam optimizer \cite{kingma2017adam} is used for the training. The learning rate follows the following rule:
\[
    \rho_k =  (0.99)^{k / 150}\, 0.001. 
\]
where the division operator stands here for the integer division, and $k$ is the train step. Thus, the learning rate is constant over 150 iterations, then decreases. It has an exponential decay with the shape of a staircase. In addition, we can reset this decay i.e. set $k=0$ at any time if we notice that the loss is reaching a plateau. The main goal of this reset strategy is to escape poor local minima of the minimization problem with a sudden large learning rate. On Figure~\ref{fig: training losses example 1} is shown a typical training and validation loss history as functions of the training step $k$ as well as the learning rate at each step. The resets enable us to make the training loss go from $1\sci{-3}$ to $5\sci{-4}$ and then to $1\sci{-4}$. 
With the loss functions weights above-mentioned, we consider $1\sci{-5}$ to be a correct plateau value to stop the training process. \answerG{The training phase lasts from 1 to 3 hours on a shared NVIDIA Tesla T4 GPU.}

\paragraph{Pre-processing} Data pre-processing is required to optimize the learning process. Here we use usual standardization techniques.

\FloatBarrier

\subsection{Numerical complexity}

\answerG{
Here, we briefly compare the computational gain in using the reduced models in the online phase. For the original system, the main cost comes from evaluating the $N$ components of the Hamiltonian gradient. If we denote by $\alpha$ the evaluation complexity for one component, the computational cost is therefore about $O(N\alpha)$. When using the reduced PSD model, an additional cost arises from the linear encoding-decoding operations and the computational cost is equal to $O(N\alpha + NK)$. The use of the DEIM-PSD method, as presented in the appendix~\ref{sec:DEIM-PSD}, allows us to rely only on $m$ components of the gradient of the Hamiltonian and the computation time is therefore about $O(m\alpha+mK)$, which no longer depends on the $N$ dimension. On the other hand, the reduced HNN model relies on the evaluation of the gradient of the neural network Hamiltonian, whose evaluation complexity is, to a first approximation, equivalent to a direct evaluation. Thus, if we denote by $n^{k}$ the width (i.e. the number of neurons) of the $k$-th layer and only count the linear operations between the layers, the total complexity of the evaluation is about $O(\sum_{k=1}^L n^{(k - 1)}n^{(k)})$. Therefore, if the width is of the order of $K^2$, the complexity of the evaluation is of the order of $O(K^4)$. Depending on the values of $\alpha$ and $m$, the reduced model AE-HNN can be competitive. It should also be noted that an additional advantage of the AE-HNN reduced model is that it can naturally be evaluated for a batch of parameters in parallel.
}

\section{Numerical results}\label{section: results}

This section is devoted to the numerical results obtained with the proposed Hamiltonian reduction method. We consider one-dimensional discretizations of three wave-type equations: the linear wave equation, the non-linear wave equation and the shallow water equation. 

\subsection{Wave equations}\label{section: Wave equations}

We introduce a parameterized one-dimensional wave equation:
\begin{align}
    \begin{cases}
        \partial_{tt} u(x, t;\mu)  - \mu_a\, \partial_x\left( w'\left( \partial_x u(x, t;\mu) , \mu_b \right) \right) + g'(u(x, t;\mu), \mu_c) = 0, &\text{ in } [0,1] \times (0,T], \\
        u(x, 0;\mu) = u_{\init}(x;\mu), &\text{ in } [0,1],\\
        \partial_t u(x, 0;\mu) = v_{\init}(x;\mu), &\text{ in } [0,1],
    \end{cases}\label{eq:non linear wave 2nd order}
\end{align}
complemented with periodic boundary conditions. The solution $u(x,t;\mu)$ represents the vertical displacement of a string over the interval $[0,1]$. The model depends on two given functions $w, g : \mathbb{R} \to \mathbb{R}$ and three parameters: $\mu=\left(\mu_a, \mu_b, \mu_c\right)^T \in \Xi \subset \R^3_+$. 

Defining the vertical displacement velocity $v(x,t;\mu) = \partial_t u(x, t;\mu)$,  Equation~\eqref{eq:non linear wave 2nd order} can be reformulated as a first order in time system:
\begin{align*}
    \begin{cases}
        \partial_t u(x, t;\mu) - v(x, t;\mu) = 0,&\text{ in } [0,1] \times (0,T], \\
        \partial_t v(x, t;\mu) + \mu_a   \partial_x\left( w'\left( \partial_x u(x, t;\mu) , \mu_b \right) \right)  + g'(u(x, t;\mu), \mu_c) = 0, &\text{ in } [0,1] \times (0,T] \\
        u(x, 0;\mu) = u_{\init}(x;\mu), &\text{ in } [0,1],\\ 
        v(x, 0;\mu) = v_{\init}(x;\mu), &\text{ in } [0,1].
    \end{cases}
\end{align*}
Then, we consider a spatial finite difference discretization of this system. Introducing a uniform mesh of the interval $[0,1]$, with $N$ cells, space step $\Delta x = 1/(N-1)$ and nodes $x_i = i \Delta x$, for $i \in \llbracket0,N-1\rrbracket$, the approximate solution $(\mathbf{u},\mathbf{v}) = (u_i(t;\mu),v_i(t;\mu))_{i \in \llbracket0,N-1\rrbracket} \in \R^{N}\times \R^N$ satisfies the following system:
\begin{align}
    \begin{cases}
        \displaystyle\frac{d}{dt} u_i (t;\mu) = v_i(t;\mu),&\text{ in } (0,T] \\[4pt]
        \displaystyle\frac{d}{dt} v_i (t;\mu) = - \frac{\mu_a}{\Delta x} \left( w'\left( \frac{u_{i+1} - u_{i}}{\Delta x} , \mu_b \right) -  w'\left( \frac{u_{i} - u_{i-1}}{\Delta x} , \mu_b \right) \right) +
        g'(u_i, \mu_c), &\text{ in } (0,T], \\
        u_i(0;\mu) = u_{\init}(x_i;\mu),& \\
        v_i(0;\mu) = v_{\init}(x_i;\mu).&
    \end{cases}\label{eq:non linear wave discrete}
\end{align}
These equations actually form a Hamiltonian system of size $2N$, with a separable Hamiltonian function ($(\mathbf{u},\mathbf{v})$ stands for variables $(\mathbf{q},\mathbf{p})$) given by
\begin{equation}\label{eq:hamiltonian function wave}
    \hami(\mathbf{u},\mathbf{v};\mu) = \hami^1(\mathbf{u};\mu) + \hami^2(\mathbf{v};\mu),
\end{equation}
with
\begin{align*}
    \hami^1(\mathbf{u};\mu)&=\Delta x \sum_{i=0}^{N-1} \left( \mu_a w \left(\frac{u_{i+1}-u_i}{\Delta x},\mu_b \right)+ \mu_a w \left(\frac{u_{i}-u_{i-1}}{\Delta x},\mu_b\right) + g(u_i, \mu_c) \right), \\
    \hami^2(\mathbf{v};\mu) &=\frac{1}{2} \Delta x \sum_{i=0}^{N-1} v_i^2.
\end{align*}
In the following, we test the AE-HNN method for two choices of functions $w$ and $g$. The same hyper-parameters are used in both test-cases. They are summarized in Table~\ref{tab:hyperparameters}.

\begin{table}[htb!]
    \centering
    \resizebox{\columnwidth}{!}{%
    \begin{tabular}{cp{4cm}ccc}
    \toprule
    && wave equations & shallow water 1D & shallow water 2D \\
    \midrule
    AE   & type &  split & bichannel & \answerG{bichannel} \\
         & nb of convolution blocks (encoder)     & $4$ & $4$ & \answerG{$4$} \\
         & dense layers (encoder) &  $[256,128,64,32]$ & $[256,128,64,32]$ & \answerG{$[512,256,128,64,32]$} \\
         &  activation functions & ELU & swish & \answerG{ELU} \\
    \midrule
    HNN & \answerG{type} & \answerG{split} & \answerG{standard} &  \answerG{standard} \\
     & dense layers & $[24,12,12,12,6]$& $[40,20,20,20,10]$ & \answerG{$[96, 48, 48, 48, 24]$} \\
        &  activation functions & $\tanh$ & swish & \answerG{$\tanh$} \\
     \midrule
     watch duration & $s$ & $16$ & $48$ & $16$ \\
    \bottomrule
    \end{tabular}}%
    \caption{Hyper-parameters. Activation functions are used except for the last layer of the neural networks. ELU refers to the function elu$(x) = x 1_{x>0} + (e^x-1) 1_{x<0}$ and swish to the function swish$(x) = x/(1+e^{-x})$. For the auto-encoder (AE), the number of convolution blocks and the sizes of the hidden of layers are those of the encoder. The decoder is constructed in a mirror way.}
    \label{tab:hyperparameters}
\end{table}

\subsubsection{Reduction of the linear wave equation}\label{section: Linear case}

\begin{figure}[htb!]
\centering
\includegraphics[width=\textwidth]{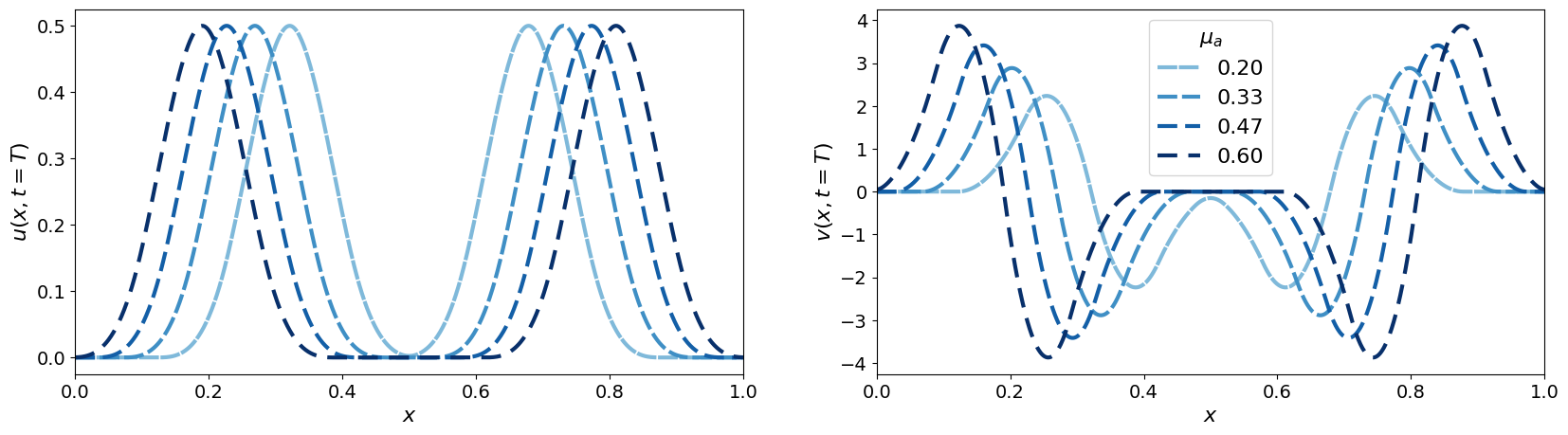}
\caption{(Linear wave) Solution $(\mathbf{u},\mathbf{v})$ at final time $T=0.4$ for various parameters $\mu=\mu_a$ }\label{fig: CI et CF linear wave}
\end{figure}

\begin{figure}
    \begin{subfigure}{\textwidth}
        \centering
        \includegraphics[width=1.\textwidth]{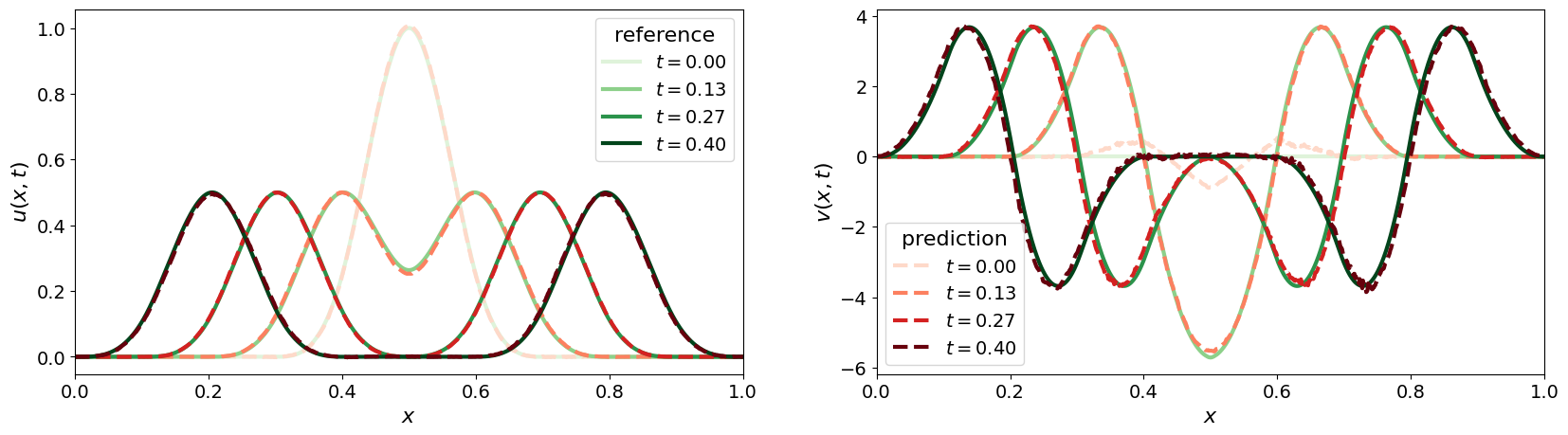}
        \caption{AE-HNN, $K=1$}
        \label{fig: linear wave figure AE HNN}
    \end{subfigure}
    \begin{subfigure}{\textwidth}
        \centering
        \includegraphics[width=1.\textwidth]{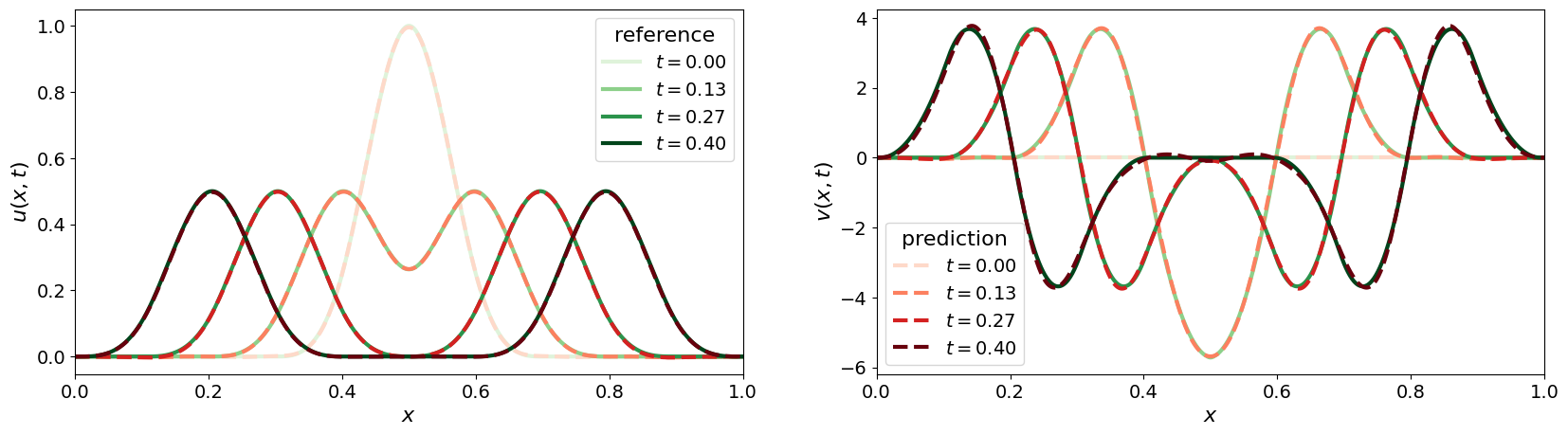}
        \caption{PSD, $K=6$}
        \label{fig: linear wave figure PSD}
    \end{subfigure}
    \caption{(Linear wave) Solution $(\mathbf{u},\mathbf{v})$ at different times on test 3, reference solution (full lines) and prediction (dashed lines).}
    \label{fig: linear wave figure}
\end{figure}

Here we consider the linear wave equation:
\begin{align*}
    \begin{cases}
    \partial_{tt} u(x, t;\mu)  - \mu_a\, \partial_{xx} u(x, t;\mu)  = 0, &\text{ in } [0,1] \times (0,T], \\
    u(x, 0;\mu) = u_{\init}(x;\mu), &\text{ in } [0,1].
    \end{cases}
\end{align*}
corresponding to $w(x, \mu_b) = \half x^2$ and $g(x,\mu_c) = 0$. In particular, we conserve only one parameter $\mu_a \in [0.2, 0.6]$, which corresponds to the square of the wave velocity. The initial condition is taken equal to:
\begin{equation}
    u_\init(x;\mu) = h(10 |x - \half|),\quad  v_\init(x;\mu) = 0. \label{eq:init_cond} 
\end{equation}
where $h$ is the compactly supported (single bump) function:
\begin{equation}h(r) = \left(1 - \tfrac{3}{2} r^2 + \tfrac{3}{4} r^3\right) 1_{[0,1]}(r) + \tfrac{1}{4} (2-r)^3 1_{(1,2]}(r).\label{eq:init_cond2}
\end{equation}
The initial condition can be observed in Figure~\ref{fig: linear wave figure AE HNN}. We consider $N=1024$ discretization points, set $T=0.4$ and $\Delta t = 1 \sci{-4}$. In Figure~\ref{fig: CI et CF linear wave}, we observe the solution at final time $T=0.4$  as a function of $\mu_a$. Each color corresponds to a different value of $\mu_a$. As we can expect, large values of $\mu_a$ generate waves farther from the initial condition.

\FloatBarrier

To train the model, we consider numerical solutions obtained with $P = 20$ different values $\mu_a$ taken regularly spaced in the interval $[0.2,0.6]$. The validation set is obtained considering numerical solutions with $6$ different parameter values in the same interval and $128$ random  pairs $(\bfY^n, \bfY^{n+s})$ for each validation parameter. The final value of the validation loss functions are given in Table~\ref{tab:validation}.

Then we test the obtained model on $3$ test values: 
\begin{equation*}
\mu_a=0.2385 \text{ (test 1),}\quad \mu_a=0.3798\text{ (test 2),}\quad \mu_a=0.5428 \text{ (test 3).}
\end{equation*}
To evaluate our model, we compute the relative $L^2$ discrete error between the reference solution $u_{\tref}$ and the reduced model prediction $u_{\pred}$:
\[
    \text{err} = \dfrac{\sum_{n=1}^{M} \Delta t\, \left(\sum_{i=1}^{N} \Delta x \, (u_{\tref, i}^n - u_{\pred, i}^n)^2\right)}{\sum_{n=1}^{M} \Delta t\, \left(\sum_{i=1}^{N} \Delta x \, ( u_{\pred, i}^n)^2\right)},
\]
where $M$ is the total number of time steps. Table \ref{tab:linear wave errors} compares the errors of the AE-HNN method with the ones obtained with the PSD and the POD for different reduced dimensions $K$. We have colored in blue the cells of the PSD and POD methods with an error lower than the AE-HNN for $K=1$ (in yellow). 
As shown in Table \ref{tab:linear wave errors}, we succeed in capturing the dynamics with a AE-HNN model of dimension $K=1$ with a relative error equal to $1\sci{-2}$. To obtain comparable performance, $K=6$ modes are needed for the PSD and $K=10$ modes for the POD. The reduction capability of the PSD is nevertheless quite good and this is due to the linearity of the problem. Also, the POD results show that taking into account the symplecticity of the problem enables us to improve significantly the reduction. 

Figure \ref{fig: linear wave figure} shows the simulations of the reduced models for both the AE-HNN $(K=1)$ and the PSD ($K=6$) methods at different times for test 3. As expected, both methods provide good results. 

\begin{table}[htb!]
    \centering
    \begin{tabular}{ccccc}
       \toprule
       & $\mathcal{L}^s_{\model}$  
       & $\mathcal{L}_{\AE}$ 
       & $\mathcal{L}^s_{\modelr}$
       & $\mathcal{L}^s_{\stab}$ \\
       \midrule
        linear wave 
        & $8.25\sci{-4}$ 
        & $8.19\sci{-4}$ 
        & $2.52\sci{-7}$
        & $1.54\sci{-5}$\\
        non-linear wave 
        &  $3.47\sci{-5}$ 
        & $3.53\sci{-5}$
        & $1.84\sci{-8}$
        & $1.98\sci{-6}$\\
        shallow water 1D
        & $6.49\sci{-5}$
        & $6.51\sci{-5}$
        & $7.70\sci{-9}$
        & $1.19\sci{-7}$ \\
        \answerG{shallow water 2D}
        & \answerG{$9.19 \sci{-5}$}
        & \answerG{$9.11 \sci{-5}$}
        & \answerG{$3.80 \sci{-8}$}
        & \answerG{$3.29 \sci{-7}$} \\
        \bottomrule
    \end{tabular}
    \caption{Values of the loss functions on validation data for the different test-cases.}
    \label{tab:validation}
\end{table}

\begin{table}[htb!]
\resizebox{\textwidth}{!}{%
\begin{tabular}{cccccccc}
\toprule
    &           & \multicolumn{2}{c}{test 1}                      & \multicolumn{2}{c}{test 2}                      & \multicolumn{2}{c}{test 3}                      \\ 
                        &  & error $u$    & error $v$    & error $u$    & error $v$    & error $u$   & error $v$    \\ \midrule
\multirow{3}{*}{AE-HNN} & $K=1$     & \cellcolor{yellow}$1.20\sci{-2}$ & \cellcolor{yellow}$6.05\sci{-2}$ & \cellcolor{yellow} $7.07\sci{-3}$ & \cellcolor{yellow} $5.84\sci{-2}$ & \cellcolor{yellow} $1.43\sci{-2}$ &\cellcolor{yellow} $8.09\sci{-2}$ \\ 
                        & $K=2$     & $6.71\sci{-3}$ & $2.24\sci{-2}$ & $1.48\sci{-2}$ & $3.99\sci{-2}$ & $9.71\sci{-3}$ & $2.30\sci{-2}$ \\ 
                        & $K=5$     & $1.67\sci{-2}$ & $2.54\sci{-2}$ & $1.48\sci{-2}$ & $2.90\sci{-2}$ & $1.79\sci{-2}$ & $3.55\sci{-2}$ \\ \midrule
\multirow{3}{*}{PSD}    & $K=4$     & $6.70\sci{-1}$ & $3.67\sci{-1}$ & $7.01\sci{-1}$ & $8.69\sci{-1}$ & $7.30\sci{-1}$ & $3.65\sci{-1}$ \\  
                        & $K=5$     & $1.28\sci{-1}$ & $1.34\sci{-1}$ & $1.60\sci{-1}$ & $1.43\sci{-1}$ & $1.91\sci{-1}$ & $1.52\sci{-1}$ \\ 
                        & $K=6$     & \cellcolor{mycolor}$4.80\sci{-3}$ & \cellcolor{mycolor}$2.11\sci{-2}$ & $5.52\sci{-2}$ & \cellcolor{mycolor}$2.14\sci{-2}$ & \cellcolor{mycolor}$5.89\sci{-3}$ & \cellcolor{mycolor}$2.23\sci{-2}$ \\ \midrule
\multirow{3}{*}{POD}    & $K=6$     & $3.07\sci{-2}$ & $1.08\sci{-1}$ & $1.30\sci{-2}$ & $2.67\sci{-2}$ & $6.00\sci{-2}$ & $1.54\sci{-1}$ \\ 
                        & $K=8$     & $1.28\sci{-2}$ & \cellcolor{mycolor}$4.89\sci{-2}$ & $1.12\sci{-2}$ & \cellcolor{mycolor}$2.66\sci{-2}$ & $3.96\sci{-2}$ & \cellcolor{mycolor}$7.05\sci{-2}$ \\ 
                        & $K=10$    & \cellcolor{mycolor}$9.78\sci{-3}$ & \cellcolor{mycolor}$4.36\sci{-2}$ & \cellcolor{mycolor}$2.30\sci{-3}$ & \cellcolor{mycolor}$9.07\sci{-3}$ & $1.58\sci{-2}$ & \cellcolor{mycolor}$5.93\sci{-2}$ \\ \bottomrule
\end{tabular}%
}
\caption{(Linear wave) Relative $L^2$ errors for different reduced dimensions $K$. Blue cells correspond to POD and PSD simulations with lower errors than the corresponding AE-HNN simulation in yellow.}
\label{tab:linear wave errors}
\end{table}

\FloatBarrier

\subsubsection{Reduction of non-linear wave equation}\label{section: results non linear wave}

Now, we consider the following non-linear wave equations:
\begin{align*}
    \begin{cases}
    \partial_{tt} u(x, t;\mu)  - \mu_a\, \partial_{xx} u(x, t;\mu) - \mu_b\, \partial_{x}  \left(\cos ( \mu_b \partial_x u(x, t;\mu)) \right) + 30 \mu_c x^2 = 0 , &\text{ in } [0,1] \times (0,T], \\
    u(x, 0;\mu) = u_{\init}(x;\mu), &\text{ in } [0,1].
    \end{cases}
\end{align*}
corresponding to $w(x, \mu_b) =  \half x^2 + \sin \left( \mu_b x \right)$ and $g(x,\mu_c) = 10 \mu_c x^3$. There are three parameters $\mu_a \in [0.2, 0.6], \mu_b \in [0.025, 0.5]$ and $\mu_c \in [0.4, 2.4]$. The initial  condition is given by Eq. \eqref{eq:init_cond}-\eqref{eq:init_cond2}. 

\begin{figure}[htb!]
\centering
\includegraphics[width=\textwidth]{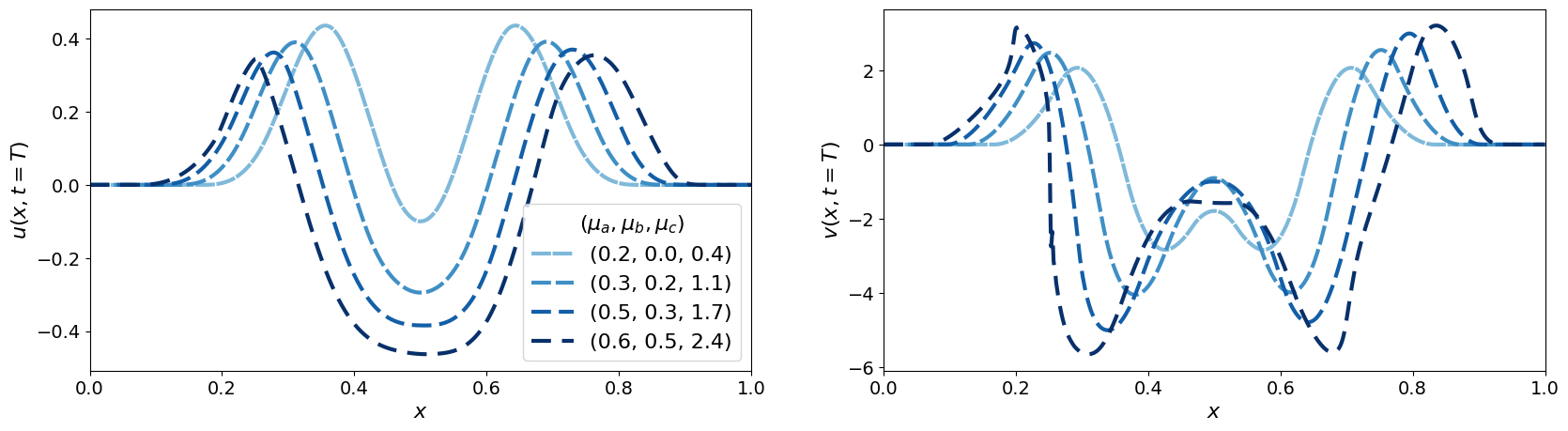}
\caption{(Non-linear wave) final condition $(\bfU,\bfV)(t=T;\mu)$ for various parameters $\mu=(\mu_a, \mu_b, \mu_c)$}\label{fig: CI et CF non linear wave}
\end{figure}

The number of discretization points is still $N=1024$ but the final time is taken equal to $T=0.3$ and the time step $\Delta t = 1\sci{-4}$. These parameters have been chosen so that the numerical simulations remain stable despite the strong non-linearities.  In Figure~\ref{fig: CI et CF non linear wave}, we observe the numerical solutions at final time for several sets of parameters. 

Training data are obtained using $P = 20$ different values $(\mu_a, \mu_b, \mu_c)$ regularly spaced in the segment $[(0.2, 0.025, 0.4), (0.6, 0.5, 2.4)]$.
The validation set is also made of $6$ different triplets define on the same domain. At the end of the training, the validation loss takes the values given in Table~\ref{tab:validation}.

The obtained model is tested on $3$ sets of parameters: 
\begin{align*}
&(\mu_a, \mu_b, \mu_c)=(0.2385, 0.088, 0.5485) \qquad \text{ (test 1),}\\
&(\mu_a, \mu_b, \mu_c)=(0.3785, 0.281, 1.354) \qquad \text{ (test 2),}\\
&(\mu_a, \mu_b, \mu_c)=(0.5528, 0.437, 2.128) \qquad \text{ (test 3),}
\end{align*}
and Table \ref{tab:non linear wave errors} presents the relative errors. A relative error of order $1\sci{-2}$ can be reached with the AE-HNN method with a reduced dimension of $K=3$ only. In comparison, the PSD and the POD require $K=15$ and $K=30$ reduced dimensions to obtain similar results. 

\begin{figure}[htb!]
    \begin{subfigure}{\textwidth}
        \centering
        \includegraphics[width=1.\textwidth]{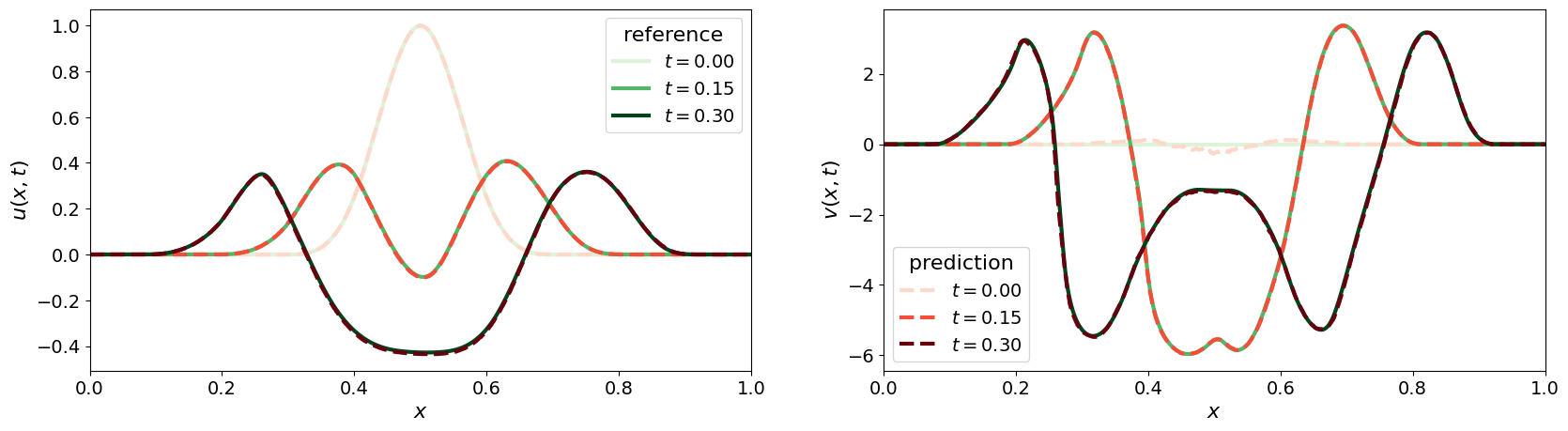}
        \caption{AE-HNN, $K=3$}
        \label{fig: non linear wave solution ae hnn}
    \end{subfigure}
    \begin{subfigure}{\textwidth}
        \centering
        \includegraphics[width=1.\textwidth]{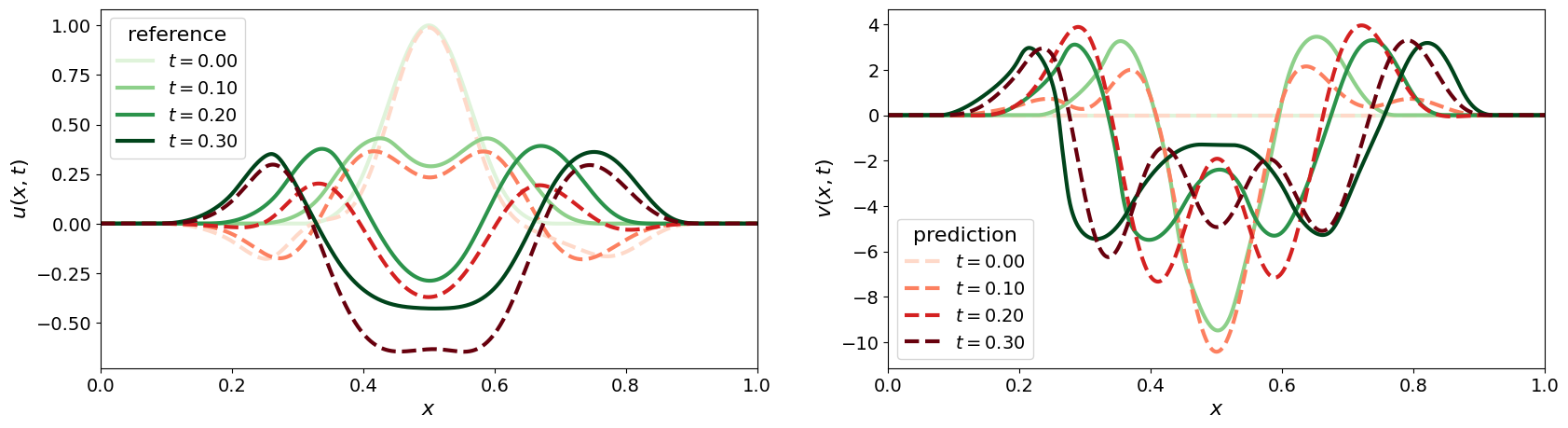}
        \caption{PSD, $K=3$}
        \label{fig: non linear wave solution psd k3}
    \end{subfigure}
    \begin{subfigure}{\textwidth}
        \centering
        \includegraphics[width=1.\textwidth]{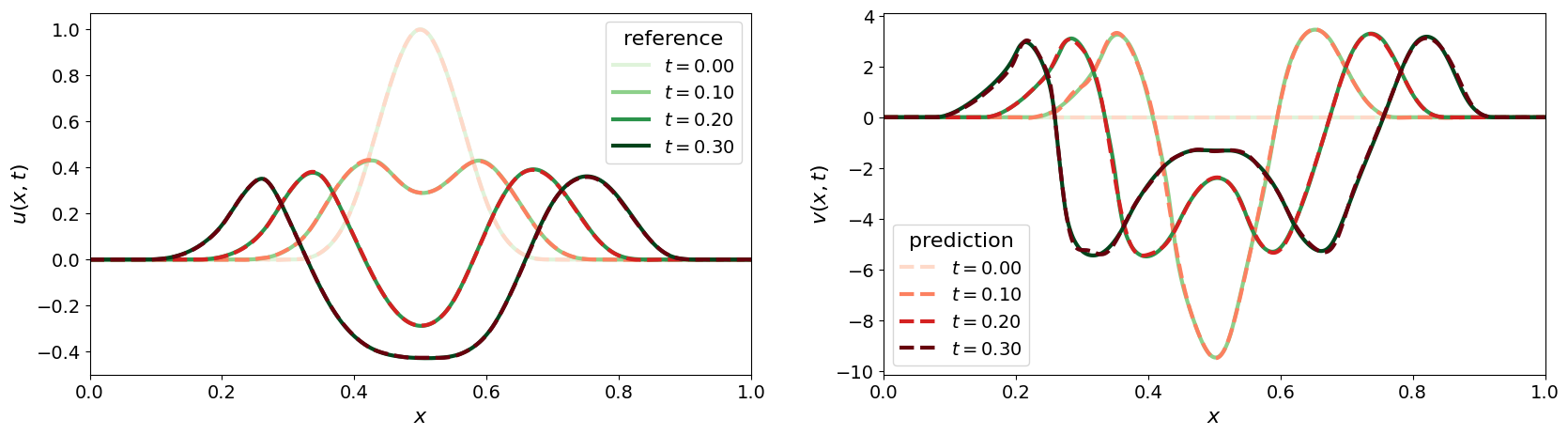}
        \caption{PSD, $K=15$}
        \label{fig: non linear wave solution psd k15}
    \end{subfigure}
    \caption{(Non-linear wave) Solutions $(\mathbf{u},\mathbf{v})$ at different times on test 3. Reference solution in solid lines and prediction in dashed lines.}
    \label{fig: non linear wave solution}
\end{figure}

Figure~\ref{fig: non linear wave solution} shows numerical solutions obtained for parameters corresponding to test 3. While the AE-HNN method with $K=3$ remains close to the reference solution (Fig. \ref{fig: non linear wave solution ae hnn}), the PSD with the same reduced dimension does not provide satisfactory results (Fig. \ref{fig: non linear wave solution psd k3}). Increasing the reduced dimension to $K=15$ for the PSD allows us to recover comparable results (Fig. \ref{fig: non linear wave solution psd k15}). 

\begin{table}[htb!]
\resizebox{\textwidth}{!}{%
\begin{tabular}{cccccccc}
\toprule
&           & \multicolumn{2}{c}{test 1}                      & \multicolumn{2}{c}{test 2}                      & \multicolumn{2}{c}{test 3}                      \\ 
                     &  & error $u$    & error $v$    & error $u$    & error $v$    & error $u$    & error $v$    \\ \midrule
AE-HNN       & $K=3$     & \cellcolor{yellow}$4.34\sci{-3}$& \cellcolor{yellow}$1.17\sci{-2}$ & \cellcolor{yellow}$5.82\sci{-2}$& \cellcolor{yellow}$8.65\sci{-2}$ &\cellcolor{yellow} $1.06\sci{-2}$& \cellcolor{yellow} $1.36\sci{-2}$ \\ \midrule
\multirow{3}{*}{PSD} & $K=3$     & $4.04\sci{-1}$& $2.55\sci{-1}$ & $4.21\sci{-1}$& $2.49\sci{-1}$ & $4.33\sci{-1}$& $2.73\sci{-1}$ \\ 
                     & $K=10$    & $3.26\sci{-2}$& $4.92\sci{-2}$ & $3.84\sci{-2}$& \cellcolor{mycolor}$4.53\sci{-2}$ & $5.63\sci{-2}$& $4.94\sci{-2}$ \\ 
                     & $K=15$    & $8.48\sci{-3}$& $1.66\sci{-2}$ & \cellcolor{mycolor}$6.45\sci{-3}$& \cellcolor{mycolor}$1.66\sci{-2}$ & \cellcolor{mycolor}$5.29\sci{-3}$& $1.87\sci{-2}$ \\ \midrule
\multirow{3}{*}{POD} & $K=3$     & $1.86\sci{-1}$& $3.16\sci{-1}$ & $5.72\sci{-2}$& $6.66\sci{-2}$ & $2.44\sci{-1}$& $2.89\sci{-1}$ \\ 
                     & $K=20$    & $2.08\sci{-2}$& $4.00\sci{-2}$ & \cellcolor{mycolor}$3.10\sci{-2}$& \cellcolor{mycolor}$7.00\sci{-2}$ & \cellcolor{mycolor}$3.66\sci{-3}$& \cellcolor{mycolor}$1.01\sci{-2}$ \\ 
                     & $K=30$    & 
                     $5.36\sci{-3}$& 
                     $2.09\sci{-2}$ & \cellcolor{mycolor}$9.54\sci{-3}$& \cellcolor{mycolor}$2.26\sci{-2}$ & \cellcolor{mycolor}$8.87\sci{-3}$& \cellcolor{mycolor}$1.93\sci{-2}$ \\ \bottomrule
\end{tabular}%
}
\caption{(Non-linear wave) Relative $L^2$ errors for different reduced dimensions $K$ and different parameters. Blue cells correspond to POD and PSD simulations with lower errors than the corresponding AE-HNN simulation in yellow.}
\label{tab:non linear wave errors}
\end{table}

\FloatBarrier


\subsubsection{Comparison with a non Hamiltonian reduction}\label{section: stability}

In this section, we want to highlight the importance of the Hamiltonian framework for the reduced dynamics. Indeed, if the Hamiltonian structure is ignored, then the reduced system can be written:
\begin{equation}
    \begin{cases}
        \displaystyle \frac{d}{dt} \bfy(t;\mu) = \overline{\mathcal F} (\bfy(t;\mu); \mu), \quad \forall t \in (0,T], \\
        \bfy(0;\mu) = \calE(\bfY_{\init}(\mu)),
    \end{cases}
\end{equation}
where $\overline{\mathcal F}:\R^K\times\Xi \to \R^K$ is the reduced vector field. Then the reduced dynamics can be learnt by designing a neural network approximation $\overline{\mathcal F}_{\theta_f}$, with a classical MLP architecture. To this aim, we can consider the following loss function:
\begin{equation}
  \mathcal{L}_{\flux}^s (\theta_e,\theta_f) = \sum_{ \bfY^{n}_{\mu} ,  \bfY^{n+s}_{\mu}  \in \, U}\left\| \bfy^{n+s}_{\mu} -  \mathcal{P}_s\left( \bfy^{n}_{\mu} ; \overline{\mathcal F}_{\theta_f} \right) \right\|^2_2,
 \label{eq:flux_loss}
\end{equation}
where the prediction operator is here a simple RK2 scheme. We discard the stability loss function $\mathcal{L}^s_{\stab}$ and keep all the other hyper-parameters identical to those of the previous section. Coupled with the auto-encoder, this reduction method will be referred to as AE-{\flux}.

We compare the two resulting method on the non-linear wave test-case considered in Section \ref{section: results non linear wave}. The HNN parameters are unchanged. For a fair comparison, the flow neural network $\mathcal F_{\theta_f}$ has the same amount of parameters: the hidden layers have $[32,24,16,16]$ units, which results in $1\,886$ parameters instead of $1\,728$ for the HNN. We train both models up to reach a loss value of about $1\sci{-5}$. The obtained validation loss functions are as follows
\begin{equation*}
    \mathcal{L}^s_{\model} = 9.46\sci{-5},\quad 
    \mathcal{L}_{\AE} = 9.88\sci{-5},\quad
    \mathcal{L}^s_{\flux} = 7.92\sci{-7}.
\end{equation*}
They are of comparable magnitude to that of the AE-HNN method in Table \ref{tab:validation}.

\begin{figure}[htb!]
    \centering
    \includegraphics[width=0.7\textwidth]{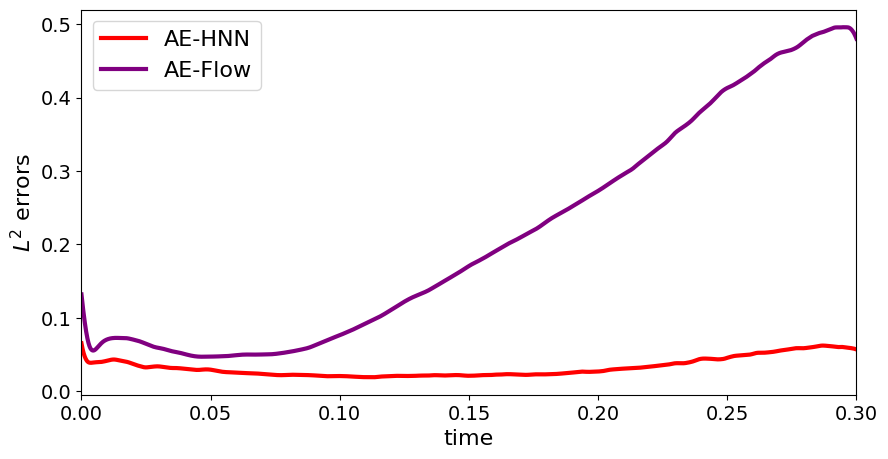}
    \caption{(Non linear wave) $L^2$ errors on test 3 as a function of time for the AE-{\flux} (purple) and the AE-HNN (red)}
    \label{fig: l2 errors vs baseline}
\end{figure}

\begin{figure}[htb!]
    \centering
    \includegraphics[width=1.\textwidth]{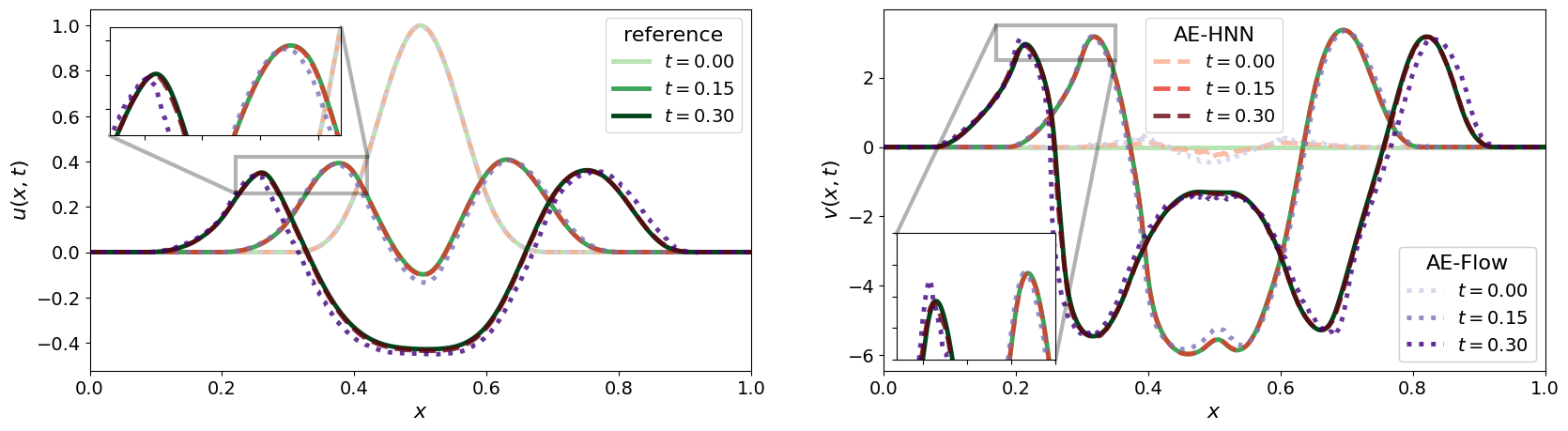}     
    \caption{(Non linear wave) $(\mathbf{u},\mathbf{v})$ for different times on test 3, reference solution (green), AE-HNN solution (red) and AE-{\flux} solution (purple)}
    \label{fig: resultat vs baseline}
\end{figure}

\begin{figure}[htb!]
    \centering
    \includegraphics[width=0.7\textwidth]{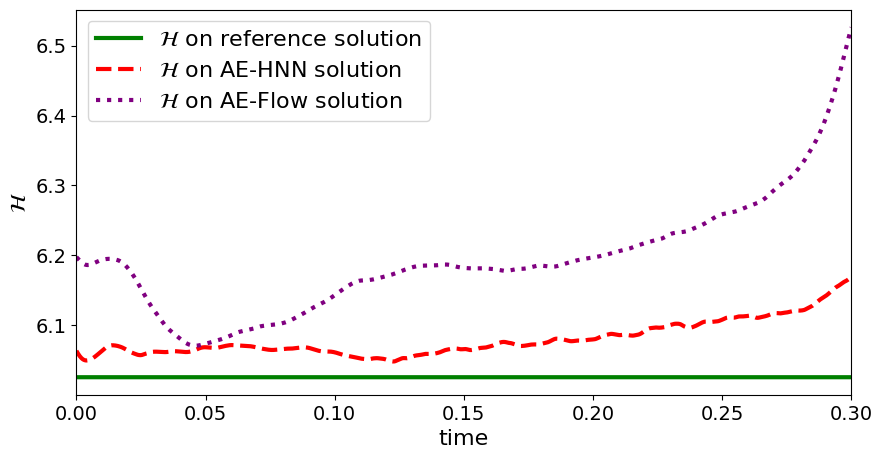} 
    \caption{(Non linear wave) Time evolution of the Hamiltonian on test 3 for the reference solution (green), the AE-HNN solution (red) and the AE-{\flux} solution (purple).}
    \label{fig: resultat vs baseline H}
\end{figure}

Relative errors are provided in Table \ref{table: relative l2 errors vs baseline}. The errors of the AE-HNN method are about $5$ times smaller than those of the AE-{\flux} method.
Figure \ref{fig: l2 errors vs baseline} shows the time evolution of the relative error for test case 3 and Figure \ref{fig: resultat vs baseline} depicts the solutions at different times. We observe that the AE-HNN solution remains close to the reference while the AE-{\flux} solution drifts from it. Furthermore, we compare in Figure \ref{fig: resultat vs baseline H} the time evolution of the Hamiltonian for the reference solution, the AE-HNN solution, and the AE-{\flux} solution on test case 3. We observe that the AE-HNN method results in a better preservation of the Hamiltonian than the AE-{\flux} solution.

\begin{table}[htb!]
\centering
\resizebox{\textwidth}{!}{%
\begin{tabular}{ccccccc}
\toprule
 & \multicolumn{2}{c}{test 1}                & \multicolumn{2}{c}{test 2}                & \multicolumn{2}{c}{test 3}                \\ 
                 & error $u$ & error $v$ & error $u$ & error $v$ & error $u$ & error $v$ \\ \midrule
AE-HNN              & $4.35\sci{-3}$ & $1.18\sci{-2}$ & $5.82\sci{-2}$ & $8.66\sci{-2}$ & $1.05\sci{-2}$ & $1.36\sci{-2}$ \\ 
AE-\flux        & $5.11\sci{-2}$ & $5.56\sci{-2}$ & $1.27\sci{-1}$ & $1.50\sci{-1}$ & $9.99\sci{-2}$ & $1.07\sci{-1}$ \\ \bottomrule
\end{tabular}%
}
\caption{(Non linear wave) Relative $L^2$ errors for the AE-{\flux} and the AE-HNN}
\label{table: relative l2 errors vs baseline}
\end{table}

\FloatBarrier

\FloatBarrier

\FloatBarrier

\subsubsection{Gain in computation time}\label{section: complexity}

In this section, we compare the actual computation time on an Intel Xeon CPU. To obtain a fair comparison, we use the same implementation of the Störmer-Verlet algorithm for all the methods. 
The only difference between the PSD and the AE-HNN methods lies in the gradient computations. For the former, we use the explicit expression of the reduced Hamiltonian
as defined in Appendix~\ref{section: PSD} and for the latter, we use the HNN backpropagation algorithm. We use 32 bits floating-point numbers and include the encoding or decoding of the data in the computation time. We consider the non-linear wave test case with the following parameters: $N=1024, K=3, T=0.4$ and $\Delta t = 1 \sci{-4}$. 

The reference solution without reduction is computed in $2\,200 \pm 14$ ms. In comparison, the AE-HNN method spends $452 \pm 20$ ms, which is five times faster. On the contrary, the PSD method spends $2\,382 \pm 12$ ms, which is a bit slower than the reference solution. This results from the evaluation of the gradient of the Hamiltonian in the non-reduced dimension. Hyper-reduction technic like Hamiltonian DEIM would accelerate the computations \cite{hesthaven2023adaptive}. However, this could also deteriorate the precision, and one should then consider larger reduced dimensions. 

\subsection{1D shallow water system}\label{sec: 1d sw}

We consider the one-dimensional shallow water system under the following formulation:
\begin{align*}
    \begin{cases}
        \partial_t \chi(x, t;\mu) + \partial_x \left( (1+\chi(x, t;\mu))\, \partial_x \phi(x, t;\mu) \right) = 0, &\text{ in } [-1,1] \times (0,T], \\
        \partial_t \phi(x, t;\mu) + \frac 1 2 (\partial_x \phi(x, t;\mu))^2 + \chi(x, t;\mu) = 0, &\text{ in } [-1,1] \times (0,T], \\
        \chi(x, 0;\mu) = \chi_{\init}(x;\mu), &\text{ in } [-1,1],\\
        \phi(x, 0;\mu) = \phi_{\init}(x;\mu), &\text{ in } [-1,1],
    \end{cases}
\end{align*}
where $\chi(x, t;\mu)$ denotes the perturbation of the equilibrium and $\phi(x, t;\mu)$ is the scalar velocity potential. Periodic boundary conditions are considered. This system admits a Hamiltonian function given by:
\[
    \mathcal H_{\text{cont}} (\chi, \phi) = \dfrac 1 2 \int_{-1}^1 \left((1 + \chi) (\partial_x \phi)^2 + \chi^2 \right) dx.
\]

Like the wave equation, we consider a spatial discretization of this model, still using a uniform mesh grid with $N$ nodes $(x_i)_{i \in \llbracket 0, N-1 \rrbracket}$ on $[-1,1]$. The approximate solution is denoted $\boldsymbol{\chi} = (\chi_i)_{i \in \llbracket 0, N-1 \rrbracket}, \boldsymbol{\phi} = (\phi_i)_{i \in \llbracket 0, N-1 \rrbracket}$. We then introduce the discrete Hamiltonian function:
\begin{equation}
    \mathcal H(\boldsymbol{\chi}, \boldsymbol{\phi}) =  \frac{\Delta x}{2} \sum_{i=1}^{N} (1 + \chi_i) \left( \frac{\phi_{i+1} - \phi_{i-1}}{\Delta x} \right)^2 + \chi_i^2,
    \label{eq: hami shallow water}
\end{equation}
and the resulting discrete shallow water equation:
\begin{equation}
    \begin{cases}
        \displaystyle \frac{d}{dt} \boldsymbol{\chi} (t;\mu) =  - D_x^2 \boldsymbol{\phi}(t;\mu) - D_x (\boldsymbol{\chi}(t;\mu) \odot D_x \boldsymbol{\phi}(t;\mu)), \\[6pt]
        \displaystyle \frac{d}{dt} \boldsymbol{\phi}(t;\mu) =  - \frac 1 2 (D_x \boldsymbol{\phi}(t;\mu)) \odot (D_x \boldsymbol{\phi}(t;\mu)) - \boldsymbol{\chi}(t;\mu),\\
        \chi_i(0;\mu) = \chi_{\init}(x_i;\mu),\\
        \phi_i(0;\mu) = \phi_{\init}(x_i;\mu), 
    \end{cases}\label{eq: FOM shallow water}
\end{equation}
where $\odot$ denotes the element-wise vector multiplication and $D_x$ the centered second-order finite difference matrix with periodic boundary condition
\[
    D_x = \frac{1}{2\Delta x}
    \begin{pmatrix}
    0 & 1 &  & -1 \\
    -1 &   & \ddots &   \\
     & \ddots &   & 1 \\
    1 &   & -1 & 0
    \end{pmatrix}.
\]
We note that Hamiltonian \eqref{eq: hami shallow water} is not separable. As a consequence, the numerical resolution of \eqref{eq: FOM shallow water} with the Störmer-Verlet scheme is implicit.Therefore, applying a reduction is all the more attractive. 

\subsubsection{Reduction of the system}\label{section: results sw}

\begin{figure}
    \begin{subfigure}{\textwidth}
        \centering
        \includegraphics[width=\textwidth]{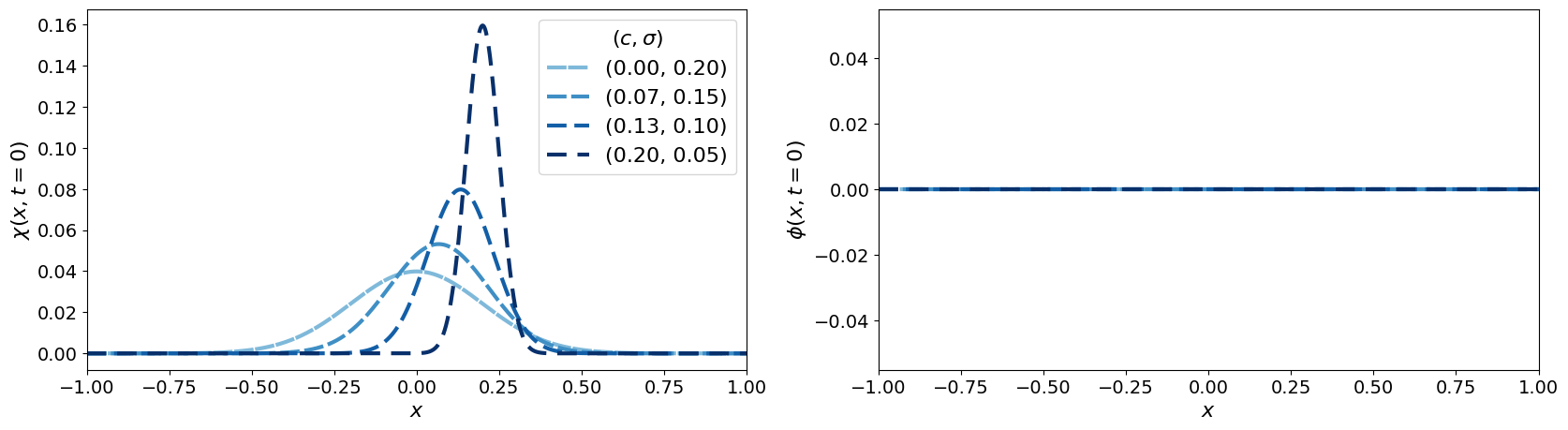}
        \caption{$t=0$}
    \end{subfigure}
    \begin{subfigure}{\textwidth}
        \centering
        \includegraphics[width=\textwidth]{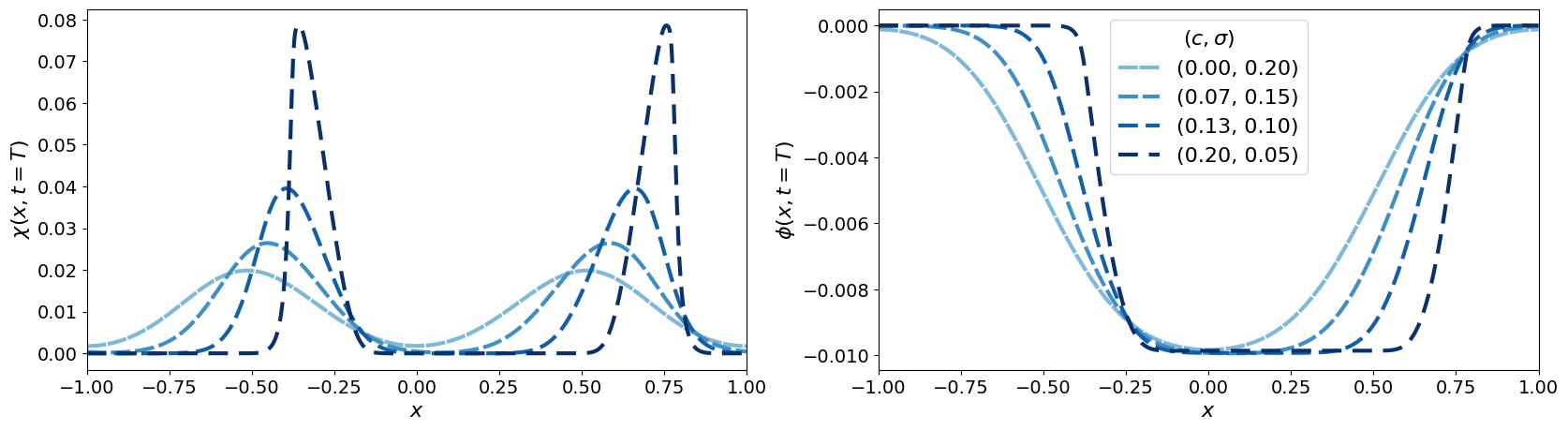}
        \caption{$t=0.5$}
    \end{subfigure}
    \caption{(Shallow water 1D) Solutions $(\boldsymbol{\chi}, \boldsymbol{\phi})$ at initial time $t=0$  and final time $t=0.5$ for various parameters $(c,\sigma)$.}\label{fig: CI et CF shallow water}
\end{figure}

The number of discretization points is taken equal to $N=1024$, the final time equal to $T=0.5$ and the time step $\Delta t = 2 \sci{-4}$.  The initial condition is parameterized with $2$ parameters $\mu = (c, \sigma) \in [0, 0.2]\times[0.2, 0.05]$ 
\begin{equation*}
    \chi_{\init}(x;\mu) = \dfrac{0.02}{\sigma \sqrt{2\pi}} \exp\left(-\dfrac{1}{2} \left(\dfrac{x - c}{\sigma}\right)^2\right),\quad \phi_{\init}(x;\mu) = 0.
\end{equation*}
For the training data, we use $20$ different solution parameters $(c, \sigma)$: we take regularly spaced values in the segment $[(0,0.2),(0.2,0.05)]$.  Let us observe the influence of the parameters on the solutions at initial and final time on Figure \ref{fig: CI et CF shallow water}. Non-linear patterns appear for small values of the standard deviation $\sigma$. 

The set of hyper-parameters is almost identical to the previous test cases except for the HNN architecture, activation functions and the watch duration, which here is taken equal to $48$. They are gathered in Table~\ref{tab:hyperparameters}.
As in the previous sections, we inspect the validation loss functions and obtain value given in Table~\ref{tab:validation}.

We choose three different sets of parameters to test the AE-HNN method:
\begin{align*}
(c, \sigma)&=(0.105, 0.11), &(\text{test 1})\\
(c, \sigma)&=(0.195, 0.053), &(\text{test 2})\\
(c, \sigma)&=(0.21, 0.045). &(\text{test 3})
\end{align*} 
As in the previous test cases, we compute the relative errors of the AE-HNN method with respect to the reference solution (see Table \ref{tab:shallow water errors}) and compare them to the results obtained with the PSD and POD methods for different values of $K$.  

The AE-HNN method achieves a mean error of about $3\sci{-2}$ with a reduced dimension of $K=4$ only. To achieve a similar performance, the PSD requires a reduced dimension of $K=24$ and the POD need a value larger than $K=32$. Table \ref{tab:shallow water errors} also shows that the AE-HNN method has a different behavior with respect to $K$ than the POD and PSD methods. For the latter, increasing $K$ improves accuracy (at the expense of computation time). With the AE-HNN method, increasing $K$ improves the encoder-decoder performance but  makes it more difficult to learn the reduced dynamics for the HNN neural network. Therefore, the HNN performance requires a balance between an adequate compression and a low-dimensional reduced model.

\begin{figure}
    \begin{subfigure}{\textwidth}
        \centering
        \includegraphics[width=1.\textwidth]{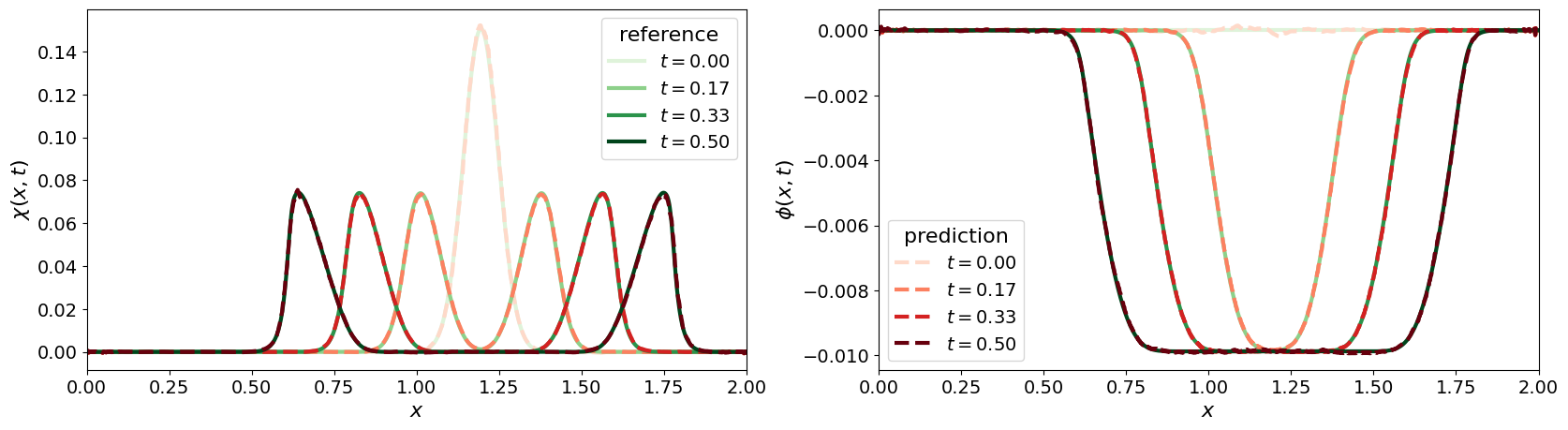}
        \caption{AE-HNN, $K=4$}
        \label{fig: resultat shallow water AE HNN}
    \end{subfigure}
    \begin{subfigure}{\textwidth}
        \centering
        \includegraphics[width=1.\textwidth]{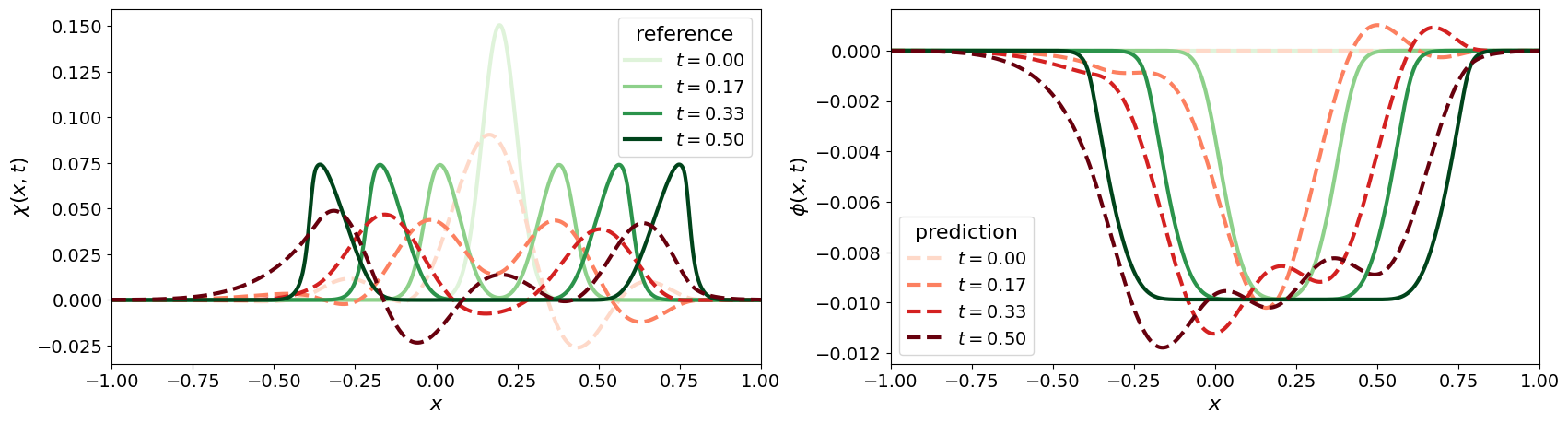}
        \caption{PSD, $K=4$}
        \label{fig: resultat shallow water PSD k4}
    \end{subfigure}
    \begin{subfigure}{\textwidth}
        \centering
        \includegraphics[width=1.\textwidth]{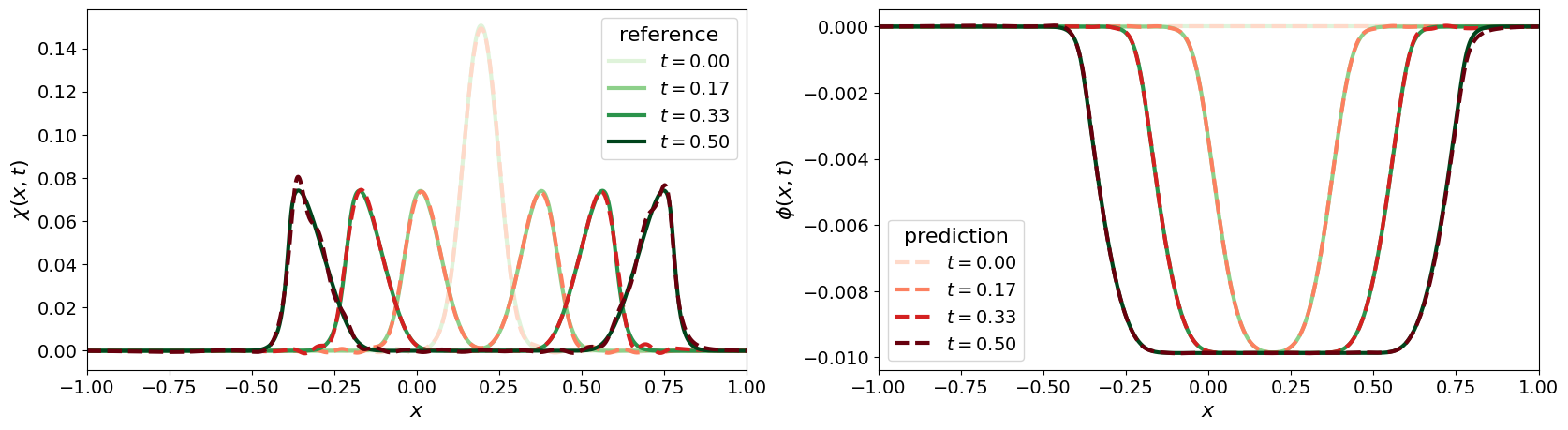}
        \caption{PSD, $K=24$}
        \label{fig: resultat shallow water PSD k24}
    \end{subfigure}
    \caption{(Shallow water 1D) $(\boldsymbol{\chi}, \boldsymbol{\phi})$ as a at different times on test 2 and $K=8$, reference solution (full lines) and prediction (dashed lines).}
    \label{fig: resultat shallow water}
\end{figure}

Then, we compare the solutions obtained with the AE-HNN method and the PSD in Figure~\ref{fig: resultat shallow water}. For a given reduced dimension $K=4$, the AE-HNN method solution remains close to the reference as expected (Fig. \ref{fig: resultat shallow water AE HNN}) while the PSD solution oscillates and stays far from the reference solution, even for the initial condition (Fig. \ref{fig: resultat shallow water PSD k4}). When increasing the dimension of the reduced model to $K=24$ for the PSD method (Fig.~\ref{fig: resultat shallow water PSD k24}), we recover similar results as the ones obtained with the AE-HNN. 


\FloatBarrier

\begin{table}[htb!]
\resizebox{\textwidth}{!}{%
\begin{tabular}{cccccccc}
\toprule
&& \multicolumn{2}{c}{test 1}                      & \multicolumn{2}{c}{test 2}                      & \multicolumn{2}{c}{test 3}                      \\ 
    &  & error $\chi$ & error $\phi$ & error $\chi$ & error $\phi$ & error $\chi$ & error $\phi$ \\\midrule
\multirow{3}{*}{AE-HNN} & $K=4$ & \cellcolor{yellow}$5.86\sci{-2}$ & \cellcolor{yellow}$2.89\sci{-2}$ & 
\cellcolor{yellow}$1.34\sci{-2}$ & 
\cellcolor{yellow}$5.00\sci{-3}$ & 
\cellcolor{yellow}$8.12\sci{-2}$ & 
\cellcolor{yellow}$2.27\sci{-2}$ \\ 
                        & $K=6$     & $8.46\sci{-2}$ & $5.17\sci{-2}$ & $2.07\sci{-2}$ & $7.77\sci{-3}$ & $9.70\sci{-2}$ & $2.80\sci{-2}$ \\ 
                        & $K=8$     & $6.26\sci{-2}$ & $3.54\sci{-2}$ & $2.21\sci{-2}$ & $1.20\sci{-2}$ & $1.45\sci{-1}$ & $4.99\sci{-2}$ \\ \midrule
\multirow{3}{*}{PSD}    & $K=10$    & $7.11\sci{-2}$ & \cellcolor{mycolor}$1.71\sci{-2}$ & $1.64\sci{-1}$ & $2.74\sci{-2}$ & $3.08\sci{-1}$ & $5.89\sci{-2}$ \\ 
                        & $K=14$    & \cellcolor{mycolor}$2.67\sci{-2}$ & \cellcolor{mycolor}$5.36\sci{-3}$ & $7.98\sci{-2}$ & $1.00\sci{-2}$ & $2.08\sci{-1}$ & $3.17\sci{-2}$ \\  
                        & $K=24$    & \cellcolor{mycolor}$1.19\sci{-2}$ & \cellcolor{mycolor}$3.43\sci{-3}$ & $2.69\sci{-2}$ & \cellcolor{mycolor}$4.20\sci{-3}$ & $9.56\sci{-2}$ & \cellcolor{mycolor}$1.06\sci{-2}$ \\ \midrule
\multirow{3}{*}{POD}    & $K=14$    & $1.13\sci{-1}$ & $4.14\sci{-2}$ & $1.22\sci{-1}$ & $4.10\sci{-2}$ & $6.98\sci{-1}$ & $2.99\sci{-1}$ \\ 
                        & $K=24$    & \cellcolor{mycolor}$4.22\sci{-2}$ & \cellcolor{mycolor}$9.91\sci{-3}$ & $3.31\sci{-2}$ & $7.81\sci{-3}$ & $3.10\sci{-1}$ & $7.96\sci{-2}$ \\  
                        & $K=32$    & \cellcolor{mycolor}$8.70\sci{-3}$ & \cellcolor{mycolor}$1.67\sci{-3}$ & \cellcolor{mycolor}$1.32\sci{-2}$ & \cellcolor{mycolor}$2.31\sci{-3}$ & $1.35\sci{-1}$ & 
                        $2.43\sci{-2}$ \\ \bottomrule
\end{tabular}%
}
\caption{(Shallow water 1D) Relative $L^2$ errors for different reduced dimensions $K$. Blue cells correspond to POD and PSD simulations with lower errors than the corresponding AE-HNN simulation in yellow.}
\label{tab:shallow water errors}
\end{table}

\FloatBarrier

\FloatBarrier

\subsubsection{Comparison with a non Hamiltonian reduction}\label{section: stability shallow water}

\begin{figure}[htb!]
    \centering
    \includegraphics[width=0.7\textwidth]{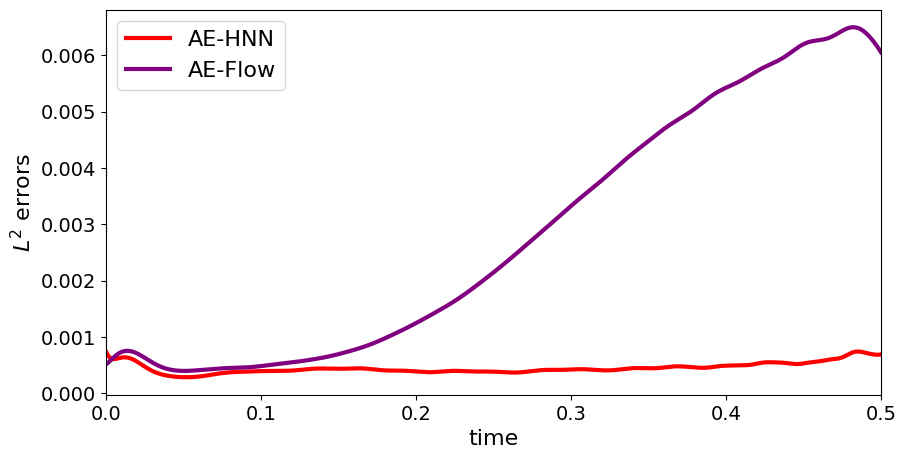}
    \caption{(Shallow water 1D) $L^2$ errors on test 2 as a function of time for the AE-{\flux} (purple) and the AE-HNN (red)}
    \label{fig: l2 errors vs baseline shallow water}
\end{figure}

\begin{figure}[htb!]
    \centering
    \includegraphics[width=1.\textwidth]{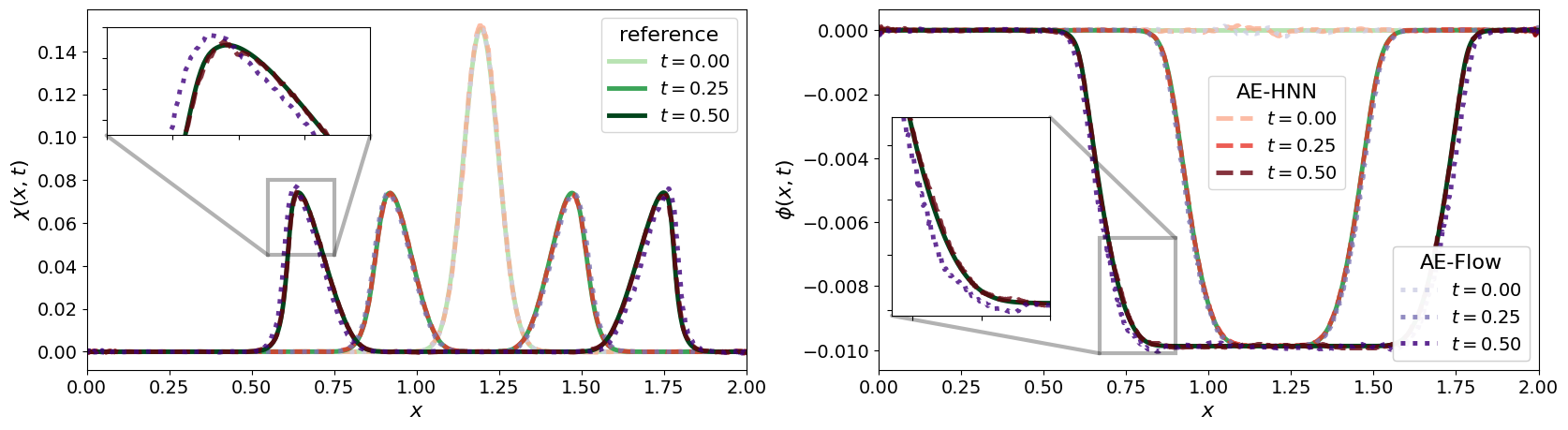}     
    \caption{(Shallow water 1D) $(\boldsymbol{\chi}, \boldsymbol{\phi})$ for different times on test 2, reference solution (green), AE-HNN solution (red) and AE-{\flux} solution (purple)}
    \label{fig: resultat vs baseline shallow water}
\end{figure}

In this section, we perform a study similar to that performed in Section~\ref{section: stability} to show the importance of having a Hamiltonian AE-HNN reduction instead of a classical AE-Flow reduced model. We consider the same test case as in Section~\ref{section: results sw} with the same hyper-parameters described in Table~\ref{tab:hyperparameters}. We stop both AE-HNN and AE-{\flux} training after reaching a validation loss value of $3\sci{-5}$. The obtained validation loss functions for the AE-{\flux} are as follows
\begin{equation*}
    \mathcal{L}^s_{\model} = 7.02\sci{-5},\quad 
    \mathcal{L}_{\AE} = 6.97\sci{-5},\quad
    \mathcal{L}^s_{\flux} = 2.85\sci{-7}.
\end{equation*}
They are of comparable magnitude to that of the AE-HNN method in Table \ref{tab:validation}. Relative errors are provided in Table \ref{table: relative l2 errors vs baseline shallow water} for the different test values. The AE-HNN method is about $4$ times more precise than the AE-{\flux} method. Figure \ref{fig: l2 errors vs baseline shallow water} shows the time evolution of the $L^2$ errors for test case  2 and Figure \ref{fig: resultat vs baseline shallow water} depicts the solution at different times. The AE-{\flux} drifts away from the reference solution while the AE-HNN remains close to it. Finally, Figure~\ref{fig: resultat vs baseline H shallow water} shows that the AE-HNN method preserves the Hamiltonian more effectively than the AE-{\flux} method.

\begin{table}[htb!]
\centering
\resizebox{\textwidth}{!}{%
\begin{tabular}{ccccccc}
\toprule
 & \multicolumn{2}{c}{test 1}                & \multicolumn{2}{c}{test 2}                & \multicolumn{2}{c}{test 3}                \\ 
                 & error $\chi$ & error $\phi$ & error $\chi$ & error $\phi$ & error $\chi$ & error $\phi$ \\ \midrule
AE-HNN              & $5.86\sci{-2}$ & $2.89\sci{-2}$ & $1.34\sci{-2}$ & $5.00\sci{-3}$ & $8.12\sci{-2}$ & $2.27\sci{-2}$ \\ 
AE-\flux        & $6.63\sci{2}$ & $3.66\sci{-2}$ & $1.05\sci{-1}$ & $3.61\sci{-2}$ & $1.73\sci{-1}$ & $5.12\sci{-2}$ \\ \bottomrule
\end{tabular}%
}
\caption{(Shallow water 1D) Relative $L^2$ errors for the AE-{\flux} and the AE-HNN}
\label{table: relative l2 errors vs baseline shallow water}
\end{table}

\FloatBarrier

\FloatBarrier

\begin{figure}[htb!]
    \centering
    \includegraphics[width=0.7\textwidth]{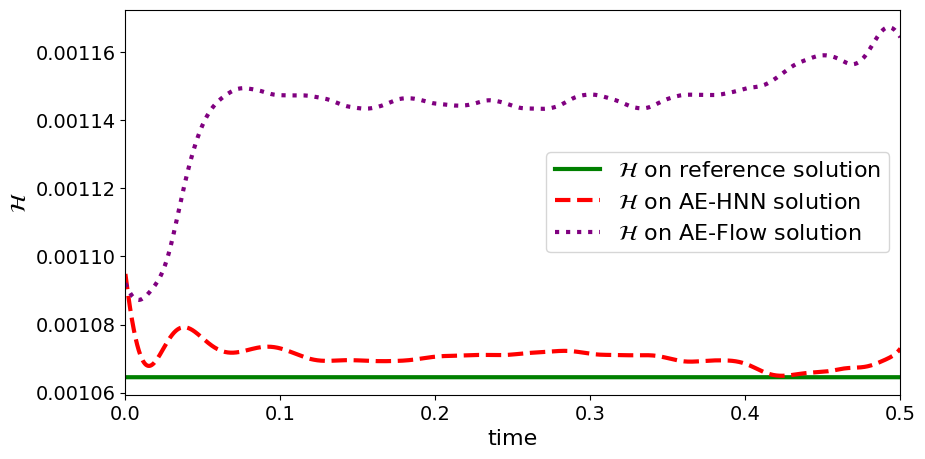} 
    \caption{(Shallow water 1D) Time evolution of the Hamiltonian on test 2 for the reference solution (green), the AE-HNN solution (red) and the AE-{\flux} solution (purple).}
    \label{fig: resultat vs baseline H shallow water}
\end{figure}
\FloatBarrier

\subsection{2D shallow water system}

\answerG{
We consider the two-dimensional shallow water system on a square  $\Omega = [-1, 1]^2$ under the following formulation:
\begin{align}
    \begin{cases}
        \partial_t \chi(\mathbf{x}, t;\mu) + \nabla \cdot \left( \left( 1 + \chi(\mathbf{x}, t;\mu) \right) \nabla \phi(\mathbf{x}, t;\mu) \right) = 0, &\text{ in }  \Omega \times (0,T], \\
        \partial_t \phi(\mathbf{x}, t;\mu) + \dfrac{1}{2} \left| \nabla \phi(\mathbf{x}, t;\mu) \right|^2 + \chi(\mathbf{x}, t;\mu) = 0, &\text{ in } \Omega \times (0,T], \\
        \chi(\mathbf{x}, 0;\mu) = \chi_{\init}(\mathbf{x};\mu), &\text{ in } \Omega, \\
        \phi(\mathbf{x}, 0;\mu) = \phi_{\init}(\mathbf{x};\mu), &\text{ in } \Omega, \\
    \end{cases}
\end{align}
where $\chi(\mathbf{x}, t;\mu)$ denotes the perturbation of the equilibrium and $\phi(\mathbf{x}, t;\mu)$ is the scalar velocity potential. Periodic boundary conditions are considered. This system admits a Hamiltonian function given by:
\[
    \mathcal H_{\text{cont}} (\chi,\phi) = \frac 1 2 \int_\Omega \left( (1+\chi) \left| \nabla \phi \right|^2 + \chi^2 \right).
\]
We consider a spatial discretization of the domain $\Omega$ with a regular mesh of $N=M^2$ cells of  size $\Delta x = 2/(M-1)$ in each direction. 
The discrete Hamiltonian function is
$$
    \mathcal H(\boldsymbol{\chi}, \boldsymbol{\phi}) = \frac 1 2
    \sum_{i,j=0}^{M-1}
    \left(
        ( 1 + \chi_{i,j})
        \left[
            \left(\dfrac{\phi_{i + 1, j} - \phi_{i - 1, j}}{2 \Delta x}\right)^2
            + \left(\dfrac{\phi_{i, j + 1} - \phi_{i, j - 1}}{2 \Delta y}\right)^2
        \right]
        + \chi_{i,j}^2
    \right)
$$
with $\chi_{i,j}(t;\mu) \approx \chi(\mathbf{x}_{i,j}, t;\mu)$ (resp. $\phi_{i,j}(t;\mu) \approx \phi(\mathbf{x}_{i,j}, t;\mu)$), with $\mathbf{x}_{i,j} = (-1,-1)+(i\Delta x, j\Delta x)$. The resulting discrete system reads
$$
    \begin{cases}
        \dfrac{d}{dt} \boldsymbol{\chi} (t;\mu) = - D_x \left( [1+\boldsymbol{\chi}(t;\mu)] \odot D_x \boldsymbol{\phi}(t;\mu) \right) - D_y \left( [1+\boldsymbol{\chi}(t;\mu)] \odot D_y \boldsymbol{\phi}(t;\mu) \right) , \\[8pt]
        \dfrac{d}{dt} \boldsymbol{\phi}(t;\mu) = - \frac 1 2 \left[ \left( D_x \boldsymbol{\phi}(t;\mu) \right)^2 + \left( D_y \boldsymbol{\phi}(t;\mu) \right)^2 \right] - \boldsymbol{\chi}(t;\mu), \\
        \chi_m(0;\mu) = \chi_{\init} (\mathbf{x}_m;\mu), \\
        \phi_m(0;\mu) = \chi_{\init} (\mathbf{x}_m;\mu),
    \end{cases}
$$
where $D_x$ and $D_y$ are respectively the centered finite difference operators along the $x$ and $y$ axis.
}

\subsubsection{Reduction of the system}\label{sec: reduction SW 2D}

\answerG{We consider $M=64$ cells per direction, The final time is set to $T=15$ and the time step to $\Delta t = 1 \sci{-3}$. In this test case, we choose to use an implicit midpoint scheme \cite{Hairer}. The initial condition, parameterized with two parameters $\mu = (\alpha, \beta) \in \left[ 0.2, 0.5 \right] \times \left[ 1, 1.7 \right]$, is chosen equal to
\[
    \chi_{\init}(\mathbf{x};\mu) = \alpha \operatorname{exp}\left(- \beta\, \mathbf{x}^T \mathbf{x}\right),\quad \phi_{\init}(\mathbf{x};\mu) = 0.
\]
For the training data, we  use $20$ different couples of parameter $(\alpha, \beta)$ regularly spaced in the segment $[\left( 0.2, 0.5 \right),\left( 1, 1.7 \right)]$. Figure~\ref{fig: shallow water 2d evol} shows the time evolution for two couples of parameters: $(\alpha, \beta) = (0.2, 1)$  and  $(\alpha, \beta) = (0.5, 1.7)$.}
\begin{figure}\ContinuedFloat
    \begin{subfigure}{\textwidth}
        \centering
        \includegraphics[width=\textwidth]{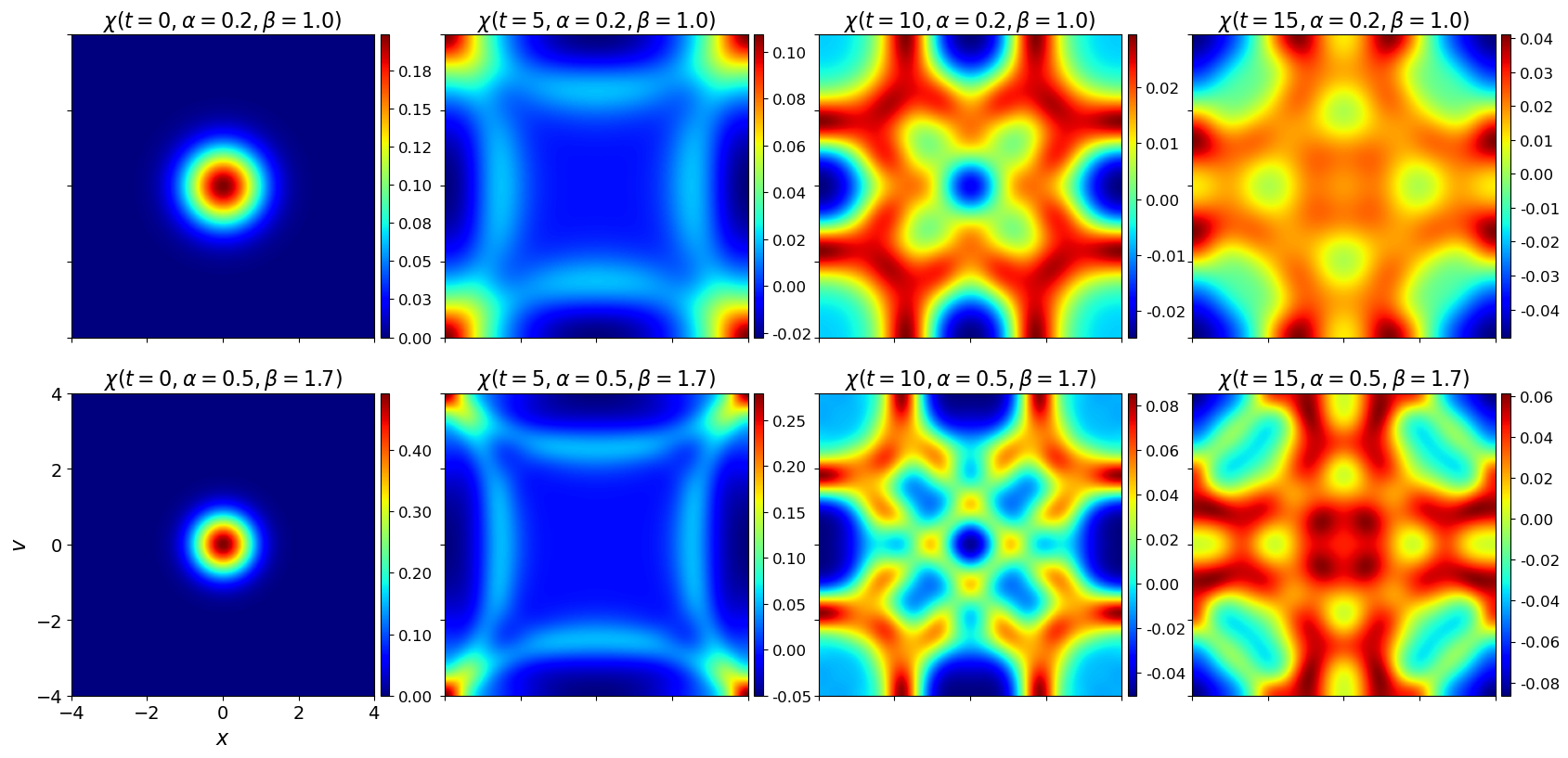}
        \caption{$\boldsymbol{\chi}(t;\mu)$}
    \end{subfigure}
    \begin{subfigure}{\textwidth}
        \centering
        \includegraphics[width=\textwidth]{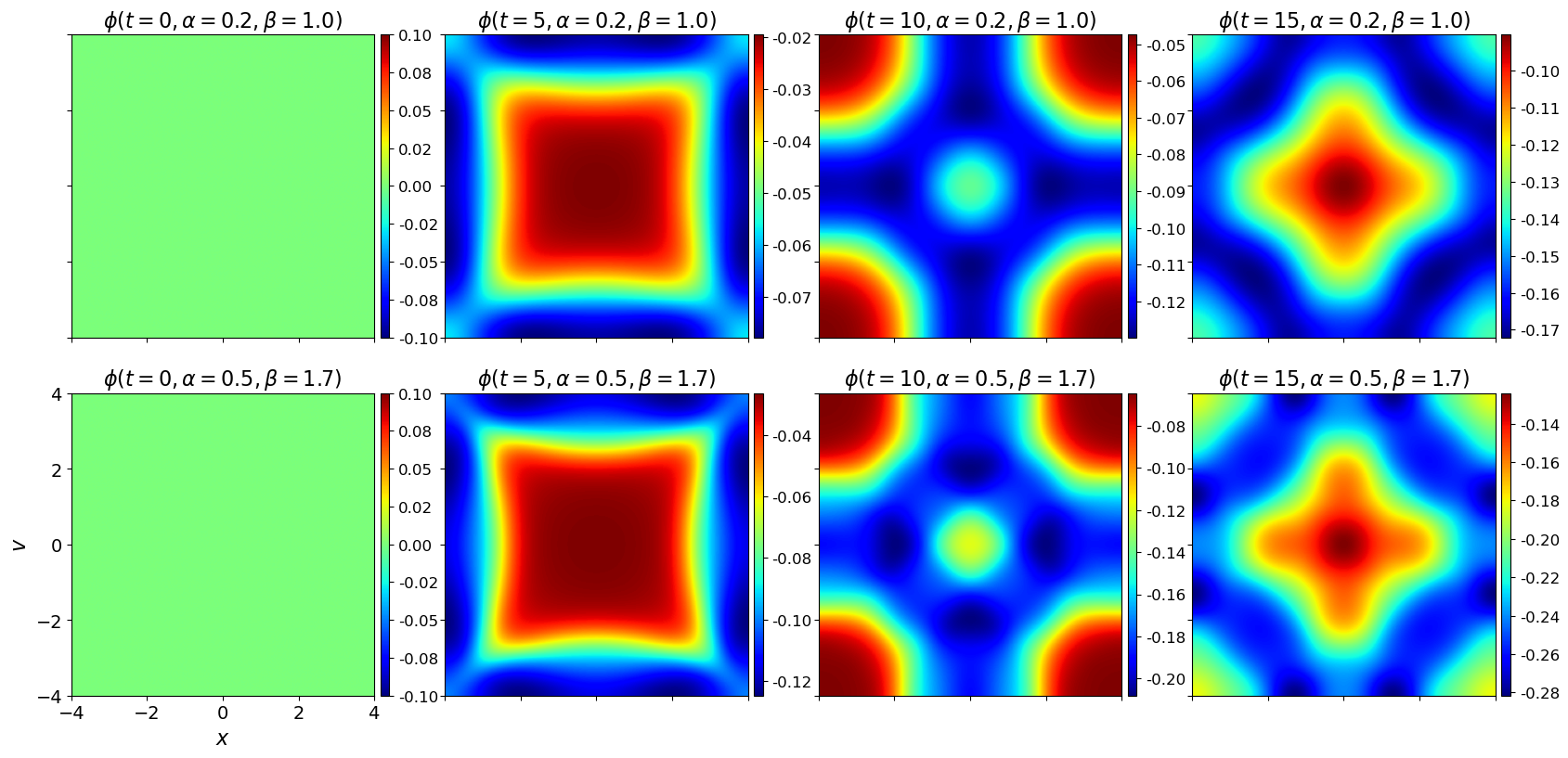}
        \caption{$\boldsymbol{\phi}(t;\mu)$}
    \end{subfigure}
    \caption{(Shallow water 2D) Solutions $(\boldsymbol{\chi}, \boldsymbol{\phi})$ at different times $t\in\{0, 5, 10, 15\}$ for various parameters $(\alpha, \beta) \in \{(0.2, 1), (0.5, 1.8)\}$.}\label{fig: shallow water 2d evol}
\end{figure}

\answerG{As data are two-dimensional, we use a 2D variant of the convolutional AE: convolutional layers are two-dimensionals, up and down-sampling are extended in 2D. All the neural network hyper-parameters are given in Table~\ref{tab:hyperparameters}. 
}

\begin{figure}[htb!]
    \centering
    \includegraphics[width=\textwidth]{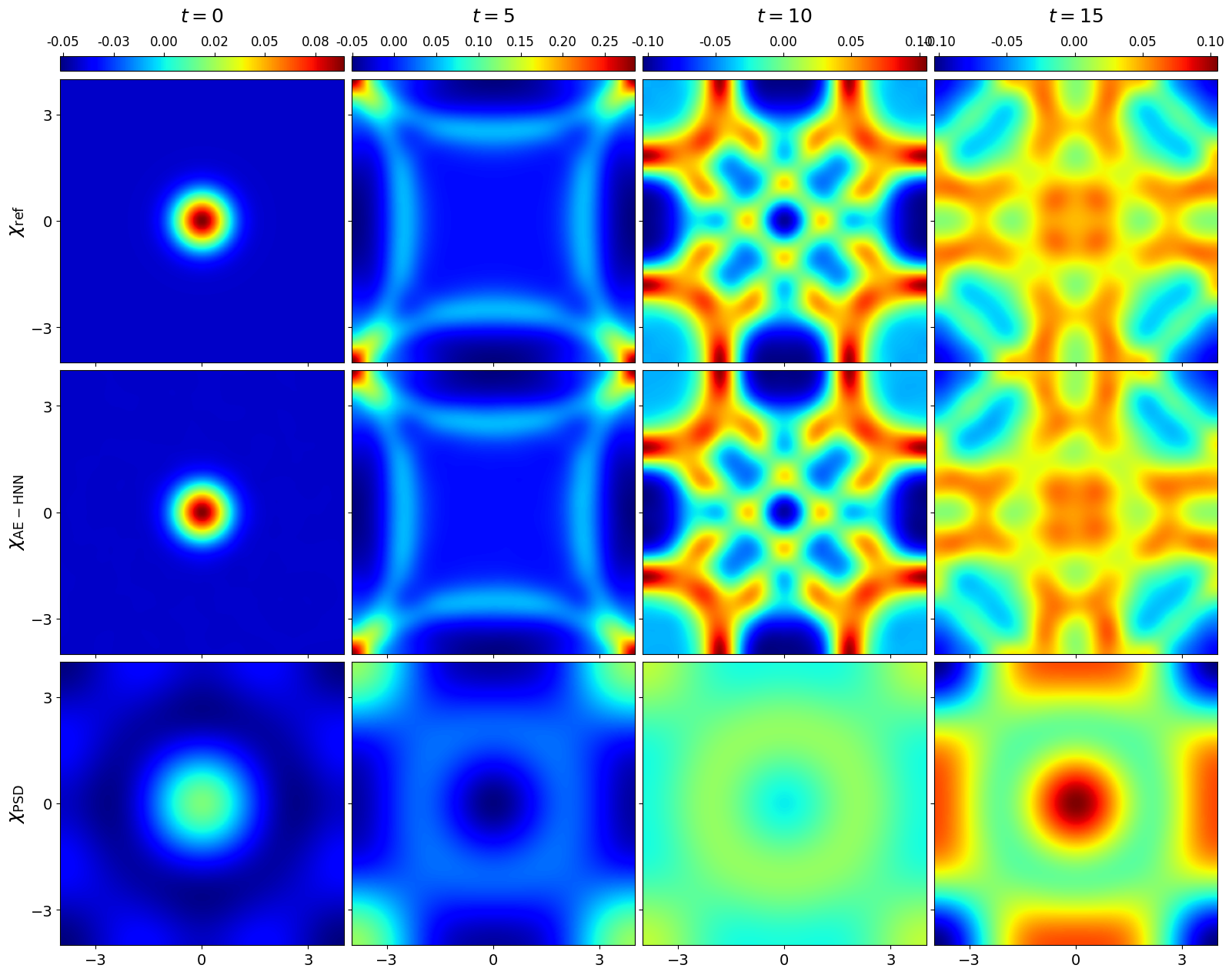}
    \caption{(Shallow water 2D) Solutions $\boldsymbol{\chi}(t;\mu)$ at different times $t\in\{0, 5, 10, 15\}$ on test 3 with $K=4$, reference solution (top line), AE-HNN solution (middle line) and PSD solution (bottom line).}\label{fig: shallow water 2d ae hnn psd}
\end{figure}

\begin{figure}[htb!]
    \centering
    \includegraphics[width=\textwidth]{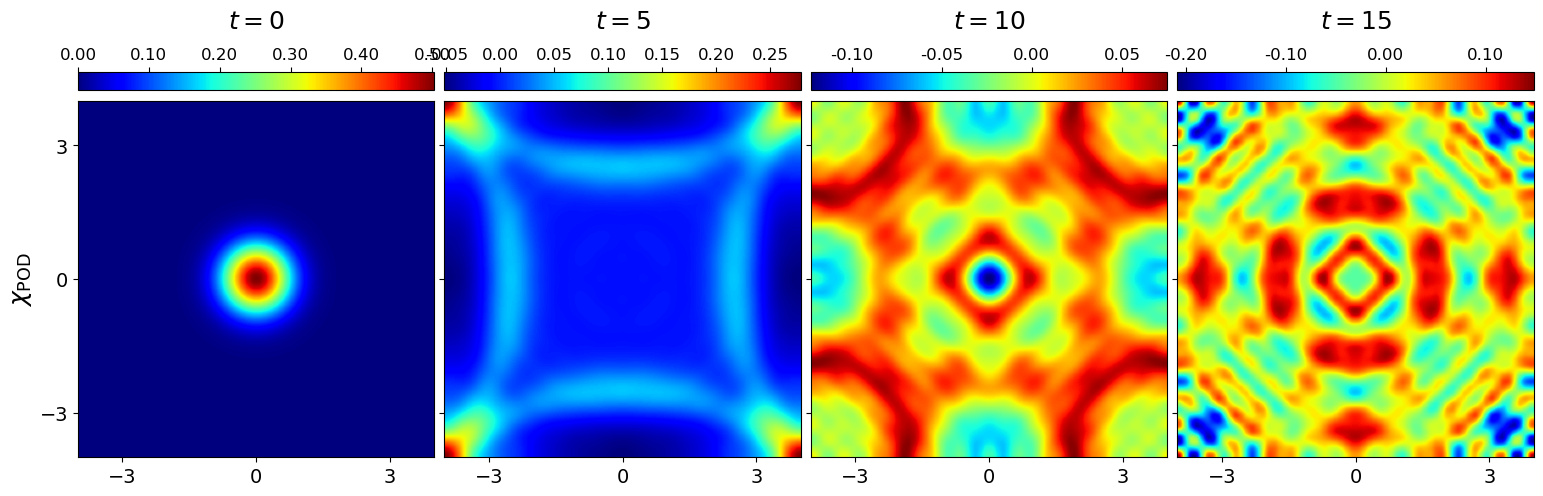}
    \caption{(Shallow water 2D) POD solution $\boldsymbol{\chi}(t;\mu)$ at different times $t\in\{0, 5, 10, 15\}$ on test 3 with $K=35$.}\label{fig: shallow water 2d pod}
\end{figure}

\answerG{
We choose three different sets of parameters to test the AE-HNN method
\begin{align*}
(\alpha, \beta)&=(0.35, 1.35), &(\text{test 1})\\
(\alpha, \beta)&=(0.41, 1.49), &(\text{test 2})\\
(\alpha, \beta)&=(0.51, 1.72). &(\text{test 3})
\end{align*} 
As in the previous test cases, we compute relative errors of the AE-HNN method with respect to the reference solution and compare them to the results obtained with the PSD and POD methods for different values of $K$ in  Table~\ref{tab:shallow water 2d errors}. Regarding the AE-HNN method, a fixed reduced dimension of $K=4$ if sufficient to obtain a precise reduced model while, for the same reduced dimension, the PSD solution is far from the reference solution as it can be observed in Fig.~\ref{fig: shallow water 2d ae hnn psd}. When increasing the reduced dimension to $K=30$, the PSD produces similar results as the ones obtained with the AE-HNN with only $K=4$. Regarding the POD reduced model, even a reduced dimension of $K=35$ does not provide the targeted precision. Indeed, the non-symplecticity of the reduced model produces some instabilities as it can be observed on test 3 in Fig.~\ref{fig: shallow water 2d pod}.
}

\begin{table}[ht]
\resizebox{\textwidth}{!}{%
\begin{tabular}{cccccccc}
\toprule
                          &      & \multicolumn{2}{c}{test 1}  & \multicolumn{2}{c}{test 2}  & \multicolumn{2}{c}{test 3}  \\
                          &      & error $\chi$ & error $\phi$ & error $\chi$ & error $\phi$ & error $\chi$ & error $\phi$ \\
\midrule
\multirow{2}{*}{AE-HNN}   & $K=4$ & \cellcolor{yellow}$4.68\sci{-2}$ & \cellcolor{yellow}$1.37\sci{-2}$  & \cellcolor{yellow}$1.73\sci{-2}$  & \cellcolor{yellow}$5.03\sci{-3}$ & \cellcolor{yellow}$3.33\sci{-2}$ & \cellcolor{yellow}$8.92\sci{-3}$  \\
                          & $K=5$ & $2.32\sci{-2}$ & $6.25\sci{-3}$  & $7.38\sci{-2}$  & $1.62\sci{-2}$ & $1.39\sci{-1}$ & $2.98\sci{-2}$  \\
\midrule
\multirow{3}{*}{PSD}      & $K=6$ & $3.48\sci{-1}$ & $3.96\sci{-2}$  & $3.82\sci{-1}$  & $4.50\sci{-2}$ & $4.33\sci{-1}$ & $5.41\sci{-2}$  \\
                          & $K=20$ & $5.05\sci{-2}$ & \cellcolor{mycolor}$3.39\sci{-3}$  & $6.52\sci{-2}$  & \cellcolor{mycolor}$4.55\sci{-3}$ & $9.49\sci{-2}$ & \cellcolor{mycolor}$7.13\sci{-3}$  \\
                          & $K=30$ & \cellcolor{mycolor}$1.37\sci{-2}$ & \cellcolor{mycolor}$9.20\sci{-4}$ & $1.87\sci{-2}$  & \cellcolor{mycolor}$1.19\sci{-3}$ & \cellcolor{mycolor}$3.09\sci{-2}$ & \cellcolor{mycolor}$2.01\sci{-3}$  \\
\midrule
\multirow{3}{*}{POD}      & $K=10$ & $4.51\sci{-1}$ & $3.69\sci{-2}$  & $4.81\sci{-1}$  & $4.47\sci{-2}$ & $5.33\sci{-1}$ & $6.02\sci{-2}$  \\
                          & $K=16$ & $1.19\sci{-1}$ & \cellcolor{mycolor}$6.31\sci{-3}$  & $4.04\sci{-1}$  & $1.76\sci{-2}$ & $1.01\sci{0}$ & $4.52\sci{-2}$  \\
                          & $K=35$ & \cellcolor{mycolor}$3.94\sci{-2}$ & \cellcolor{mycolor}$1.83\sci{-3}$  & $5.74\sci{-2}$  & \cellcolor{mycolor}$2.67\sci{-3}$ & $5.05\sci{-1}$ & $2.52\sci{-2}$  \\
\bottomrule
\end{tabular}}%
\caption{(Shallow water 2D) Relative $L^2$ errors for different reduced dimensions $K$. Blue cells correspond to POD and PSD simulations with lower errors than the corresponding AE-HNN simulation in yellow.}
\label{tab:shallow water 2d errors}
\end{table}

\FloatBarrier

\subsubsection{Comparison with a symplectic DEIM hyper-reduction}

\answerG{
While the AE-HNN would require a smaller reduced dimension than the PSD for a given targeted precision, practical efficiency still need to be compared. As the PSD still requires to come back to the original $2N$ dimension for non-linear models, we also compare the AE-HNN with the Discrete Empirical Interpolation Method (DEIM) version of the PSD as proposed in \cite{peng2015symplectic}: it relies on an interpolation of $m$ given components (among the $2N$ ones). For the sake of completeness, the method is briefly presented in Annex~\ref{sec:DEIM-PSD}.}

\begin{figure}[htb!]
    \centering
    \includegraphics[width=\textwidth]{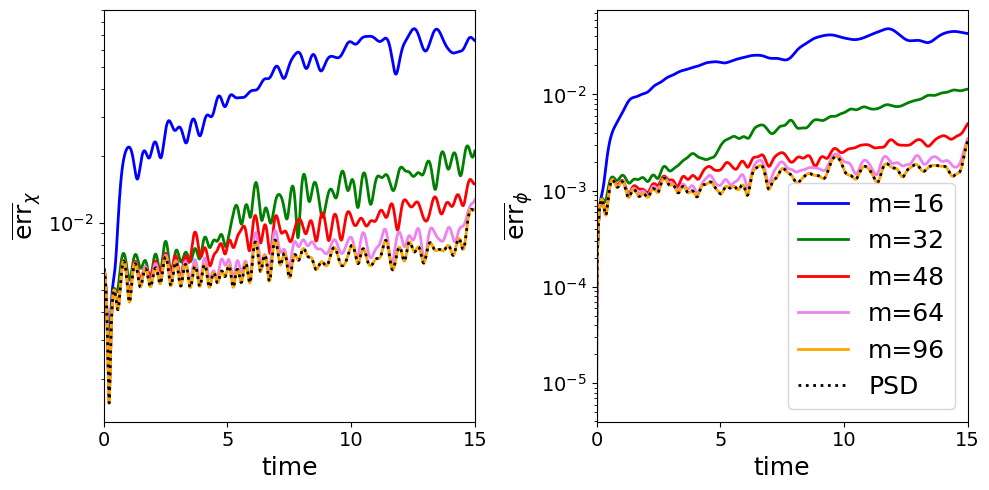}
    \caption{(Shallow water 2D) $\overline{\text{err}}$ for solutions $\boldsymbol{\chi}(t;\mu)$ (left) and $\boldsymbol{\phi}(t;\mu)$ (right) for different $m\in\{0, 5, 10, 15\}$ with $K=30$.}\label{fig: shallow water 2d mean error k30}
\end{figure}

\answerG{First, we test the precision of this DEIM approximation as a function of $m$. We consider the test-case of the previous subsection, with reduced dimension equal to $K=30$. We compute the discrete $L^2$ error on the 3 test cases:
\[
    \overline{\text{err}}_u^n = \frac{1}{3}\sum_{\mu \in \text{ test }}\left(\sum_{i=1}^{N} \Delta x \, (u_{\tref, i}^{\mu, n} - u_{\pred, i}^{\mu,n})^2\right).
\]
where $\u_{\tref, i}$ refers to the reference solution computed with the original model and $u_{\pred, i}^{\mu,n}$ the one computed with the reduced one.
We observe the time evolution of this error for different values of $m$ on Fig~\ref{fig: shallow water 2d mean error k30}. As expected, the DEIM method deteriorates the reduced model solution with respect to the sole PSD 
and increasing $m$ decreases the error.  
In practice, the value of $m$ is set so that the DEIM does not deteriorate the solution too much.}

\answerG{We then compare the computational times of the AE-HNN method with respect to the PSD algorithm with and without symplectic DEIM hyper-reduction. We set $K=4$ and $m=32$ so that the DEIM error is smaller than the PSD error. Numerical integration are performed with an implicit midpoint scheme as above-mentioned. The original model and the PSD without DEIM are solved on an Intel Xeon CPU while the AE-HNN method is executed on a NVIDIA Tesla T4 GPU. We also use the Strang splitting method \cite{Hairer} where the linear part is solved explicitly and the non-linear one is solved with numerical integration except for the AE-HNN method which cannot be split. 
The computational time equals $101.0$ s for the original model, $107.0$ s for the PSD and $57.4$ s for the DEIM-PSD.  
As expected, the PSD is not time efficient in a non-linear case due to the decompression-compression of the solution at each time step. An efficient DEIM-PSD implementation allows a satisfactory speed-up of a factor $1.76$. Regarding the AE-HNN model, it is solved in $26.5$ s, which corresponds to a speed up of $3.81$. This could be further improved. Indeed, the neural networks are executed on a GPU but the non-linear solver used in the implicit part of the scheme is executed on a CPU, which slows down the reduced order model. In practice, a SciPy \cite{2020SciPy-NMeth} implementation of the nonlinear solver \cite{article_df_sane} is used. More than $62\%$ of the non-linear solver computation time takes place on the CPU. 
}


\section{Conclusion}

We have developed a new Hamiltonian reduction method. It is based on an auto-encoder (AE) to transform initial variables into reduced variables and vice versa, and on a Hamiltonian Neural Network (HNN) to learn the Hamiltonian reduced dynamics. Using a set of coupled loss functions, we are able to learn a reduced model (AE-HNN), which has an Hamiltonian structure. It already has better reduction properties than the PSD in linear test-cases, but the gain is much larger in non-linear test-cases, as expected. Due to its Hamiltonian structure, the reduced model also shows good stability properties. \answerG{Two-dimensional test-cases show that the AE-HNN method has better computational performance than the PSD-DEIM method.}

\answerG{The question remains of how to improve the quality of the approximation. Indeed, unlike the PSD method, increasing the reduced dimension $K$ does not always provide better results.  The quality of the approximation could rather be increased by modifying the architecture of the neural networks (increasing the number of layers and the number of neurons per layer), which implies increasing the number of trainable parameters. However, up to our knowledge, there is unfortunately a lack of theoretical results about a systematic way to improve such approximations. Further studies need to be carried out.}

Obviously, the results will have to be extended to partial differential equations in three spatial dimension. As convolution layers are used in the auto-encoder, the increase of dimension should not require too large an increase in the size of the neural networks. In consequence, we expect that the computational gain will be even larger for such dynamics.  Other extensions would be to consider time-dependent reductions as in \cite{hesthaven2023adaptive} and adapt the method for the reduction of large Hamiltonian differential equations that do not have spatial structures \cite{hesthaven2023adaptive}.

\paragraph{Acknowledgements.} This research was funded in part by l’Agence Nationale de la Recherche (ANR), project
ANR-21-CE46-0014 (Milk). 


\appendix

\section{Linear reduction}\label{section: PSD}

\subsection{Proper Symplectic Decomposition}\label{sec:PSD}
Here we briefly present the Proper Symplectic Decomposition (PSD) \cite{peng2015symplectic}. In this approach, we assume that the decoder is a linear operator 
\[ \calD_{\theta_d}(\bfy) = A \bfy, \]
with $A \in M_{2N,2K}(\R)$. To be consistent with the notations, we keep the subscript $\theta_d$ to refer to the decoder parameters, which here are the coefficients of the matrix $A$.  The trial manifold $\widehat{\mathcal{M}}$ is a hyperplane. We further assume that $A$ is symplectic, i.e. satisfies $A^T J_{2N} A = J_{2K}$, and we denote by $\text{Sp}_{2N,2K}(\R)$ the set of symplectic matrices. The encoder is thus defined as the symplectic inverse of $A$: 
\[ \calE_{\theta_e}(\bfY) = A^+ \bfY, \]
with $A^+ = J_{2K}^T A^T J_{2N} \in \text{Sp}_{2K,2N}(\R)$ and where $\theta_e = \theta_d$ also denotes the coefficients of the matrix $A$. Then the projection operator onto $\widehat{\mathcal{M}}$ writes $A A^+$ and $A$ is determined by solving the minimization problem
\[
    \underset{A \in \text{Sp}_{2N,2K}(\R)}{\operatorname{min}}\  \| Y - A A^+ Y \|^2_F,
\]
where $Y$ refers to the snapshot matrix:
\begin{equation}
    Y = \begin{bmatrix}
        \bfY_1, \ldots, \bfY_p 
    \end{bmatrix} \in M_{2N, p}(\R),
    \label{eq:Y}
\end{equation}
where $(\bfY_i)_{i \in \llbracket 1,p \rrbracket}$ are $p$ values obtained  from numerical simulations of the original problem at different times and for different parameters, and $\|X\|_F = \sqrt{\sum_{ij} X_{ij}^2}$ is the Frobenius norm of $X \in M_{2N, p}(\R)$.

An approximate solution to this problem can be obtained with the cotangent lift algorithm \cite{peng2015symplectic}. First, the snapshot matrix $Y$ is reshaped into
\begin{equation*}
    \tilde Y = \begin{bmatrix}
        \bfQ_1, \ldots, \bfQ_p,
        \bfP_1, \ldots, \bfP_p
    \end{bmatrix} \in \mathcal{M}_{N,2p}(\R).
\end{equation*}
Its Singular Value Decomposition (SVD) is computed leading to  $\tilde Y =  U \Sigma V^T$, where $\Sigma$ is the diagonal matrix of the singular values and $U$ (resp. $V$) the matrix of the left (resp. right) singular vectors. Then, we define  $\Phi = U[:,:\!\! K]$ composed of the $K$-th left eigenvectors associated with the $K$-th largest singular values and we consider the following orthogonal symplectic matrix: 
\[
    A = \begin{pmatrix}
        \Phi & 0 \\ 0 & \Phi
    \end{pmatrix}.
\]
The reduced model satisfied by the reduced variable $\bfy = A^+ \bfY$ can be obtained by multiplying \eqref{eq: symplectic formulation} by $A^+$ so that
\begin{align*}
\displaystyle \frac{d}{dt} \bfy &= A^+ J_{2N} \, \nabla_{\mathbf{y}} \hami ( \bfY)  = J_{2K} A^T  \, \nabla_{\mathbf{y}} \hami ( \bfY) \approx J_{2K}  \, \nabla_{\bfy} (\hami \circ A) \left( \bfy\right),
\end{align*}
where we use the expression of $A^+$, the identity $J_{2N}^2 = - \text{Id}$ and the approximation $\bfY \approx A \bfy = A A^+ \bfY$. The reduced model involves the following reduced Hamiltonian  
\[ \overline{\hami}_{\theta_h} = \hami \circ A, \]
where $\theta_h = \theta_e = \theta_d$ still denotes the coefficients of the matrix $A$, and whose evaluation requires to compute the non-reduced quantity $A\bfy \in \R^N$. 
To avoid this evaluation in large dimension,  approximations of the reduced Hamiltonian can be constructed using the discrete empirical interpolation method (DEIM) \cite{wang2021structurepreserving,peng2015symplectic}. 

The non-symplectic version of the algorithm is called the Proper Orthogonal Decomposition (POD) \cite{papier_pod}, where the matrix $A$ is obtained using the $K$ first singular vectors of the SVD of the original snapshot matrix $Y$ defined in \eqref{eq:Y}.

\subsection{Symplectic Discrete Empirical Interpolation Method}
\label{sec:DEIM-PSD}

\answerG{
In this appendix, we shortly present the symplectic version of the Discrete Empirical Interpolation Method (DEIM) proposed in \cite{peng2015symplectic}. We suppose that the Hamiltonian function $\mathcal{H}$ can be split into a linear and a non-linear part. Its gradient then reads as
\[
    \displaystyle \nabla_\bfY \mathcal{H}(\bfY) = L \bfY + f_N(\bfY).
\]
where $L \in \mathcal{M}_{2N, 2N}(\R)$ and $f_N$ is a non-linear function. 
 As in Sec.~\ref{sec:PSD}, we apply a symplectic Galerkin projection with a linear decoder operator $A \in \mathcal{M}_{2N, 2K}(\R)$ to derive the reduced model
\begin{equation*}
    \displaystyle \frac{d}{dt} \bfy = J_{2K} A^T  \, \nabla_{\mathbf{y}} \hami ( A\bfy) =  J_{2K} (A^T L A) \bfy + J_{2K} A^T f_N(A \bfy)
\end{equation*}
While the matrix $(A^T L A)$ can be pre-computed, the non-linear part cannot be simplified and the computational cost still depends on the initial dimension $N$. The DEIM strategy consists into approximating this term by using $m$ components of $f_N$ only. The selected components are denoted $\beta_1, \ldots, \beta_m \in \llbracket1,2N\rrbracket$ and the selection matrix is $P = (\mathbf{e}_{\beta_1}, \ldots, \mathbf{e}_{\beta_m}) \in \mathcal{M}_{2N,m}(\R)$  where $\mathbf{e}_{\beta_i}$ is the $\beta_i$-th element of the canonical basis.  The approximation is then done as follows: first consider $\Psi \in \mathcal{M}_{2N,m}(\R)$, obtained after performing a SVD of samples of $f_N(\bfY)$ (as presented in Section \ref{section: PSD}) and retaining only the $m$-th largest modes. Then, the indices $\beta_1, \cdots, \beta_m$ are selected with a greedy algorithm applied on $\Psi$ and described in \cite{peng2015symplectic}. Next, we use the following approximation 
\begin{equation*}
    f_N(\bfY) \approx \Psi \left[\left(P^T \Psi\right)^{-1} P^T f_N(\bfY)\right].
\end{equation*}
This approximation is chosen so that it belongs to the range of $\Psi$ and the $m$ selected components exactly coincide with those of $f_N(\bfY)$, which can be checked after multiplication by $P^T$. The matrix $(P^T \Psi)$ is assumed to be invertible. Hence, the reduced model reads
\[
    \displaystyle \frac{d}{dt} \bfy(t) =  J_{2K}(A^T L A) \bfy + J_{2K} \left(A^T\Psi \left(P^T \Psi\right)^{-1}\right) g_m(\bfy)
\]
with $g_m(\bfy) = P^T f_N(A\bfy)$. The matrix $(A^T\Psi \left(P^T \Psi\right)^{-1})$ can be precomputed. If each component of $f_N$ depends only on a few number of components of $A\bfy$, then the evaluation complexity of $g_m$ becomes a function of $m$ instead of $N$.  }

\bibliographystyle{plain}
\bibliography{references}

\end{document}